\documentclass[11pt,a4paper]{amsart}

\usepackage{amssymb, amsmath, amsthm, amscd}

\usepackage{enumerate}

\theoremstyle{plain} \newtheorem{theorem}{Theorem}[section]
\theoremstyle{plain} \newtheorem{corollary}[theorem]{Corollary}
\theoremstyle{plain} \newtheorem{proposition}[theorem]{Proposition}
\theoremstyle{plain}\newtheorem{lemma}[theorem]{Lemma}
\theoremstyle{definition} \newtheorem{definition}[theorem]{Definition}
\theoremstyle{definition}\newtheorem{example}[theorem]{Example}
\theoremstyle{remark}\newtheorem{remark}[theorem]{Remark}
\theoremstyle{definition}\newtheorem{hypothesis}[theorem]{Hypothesis}
\theoremstyle{remark}
\theoremstyle{definition}
\theoremstyle{definition}

\def\restriction#1#2{\mathchoice
              {\setbox1\hbox{${\displaystyle #1}_{\scriptstyle #2}$}
              \restrictionaux{#1}{#2}}
              {\setbox1\hbox{${\textstyle #1}_{\scriptstyle #2}$}
              \restrictionaux{#1}{#2}}
              {\setbox1\hbox{${\scriptstyle #1}_{\scriptscriptstyle #2}$}
              \restrictionaux{#1}{#2}}
              {\setbox1\hbox{${\scriptscriptstyle #1}_{\scriptscriptstyle #2}$}
              \restrictionaux{#1}{#2}}}
\def\restrictionaux#1#2{{#1\,\smash{\vrule height .8\ht1 depth .85\dp1}}_{\,#2}}
\numberwithin{equation}{section}

\usepackage{color}
\usepackage{graphicx}
\usepackage{epsfig}

\usepackage{float}

\title[Moduli spaces for CR submersions]{Kuranishi type Moduli Spaces for proper CR submersions fibering over the circle}
\author{Laurent Meersseman}
\date{\today}
\thanks{This work was partially supported by project COMPLEXE (ANR-08-JCJC-0130-01) from the Agence Nationale de la Recherche. 
It is part of Marie Curie project DEFFOL 271141 funded by the European Community. I enjoyed the warmful atmosphere of the CRM of Bellaterra during the elaboration of this work. I would like to thank the ANR, the CRM and the European Community. This work benefited from fruitful discussions with Marcel Nicolau to whom I am especially grateful.}

\subjclass{32G07, 
32V05, 
57R30
}

\address{Laurent Meersseman\\ LAREMA\\ Universit\'{e} d'Angers\\ F-49045 Angers Cedex, France\\ laurent.meersseman@univ-angers.fr}

\begin{document}
\begin{abstract}
Kuranishi's fundamental result (1962) associates to any compact complex manifold $X_0$ a finite-dimensional analytic space which has to be thought of as a local moduli space of complex structures close to $X_0$. In this paper, we give an analogous statement for Levi-flat CR manifolds fibering properly over the circle by associating to any such $\mathcal X_0$ the loop space of a finite-dimensional analytic space which serves as a local moduli space of CR structures close to $\mathcal X_0$. We then develop in this context a Kodaira-Spencer deformation theory making clear the likenesses as well as the differences with the classical case. The article ends with applications and examples.
\end{abstract}

\maketitle

\section{Introduction.}
\label{intro}
In 1962, M. Kuranishi proved that any compact complex manifold $X_0$ has a versal (also called semi-universal) finite-dimensional analytic space of deformations $K$ (see \cite{Ku1} for the original paper, \cite{Ku2} and \cite{Ku3} for more accurate versions and simpler proofs). Roughly speaking, this means that every deformation of the complex structure of $X_0$ is encoded in a family defined over this finite-dimensional analytic space $K$ (the family is said to be complete); and that this family is minimal amongst all complete families. This fundamental result is the crowning achievement of the famous deformation theory K. Kodaira developed in collaboration with D.C. Spencer (see \cite{K-S1}, \cite{K-S2}). The Kuranishi space $K$ must be thought of as a substitute for a local moduli space of complex structures, which is known not to exist in general in the category of analytic spaces. Together with the theorems of deformations of Kodaira-Spencer (criteria of rigidity, completeness and versality in terms of the Kodaira-Spencer map), it gives a complete answer to the original Kodaira's problem \cite[p.19]{MK}: Determine all "sufficiently small" deformations of a given complex manifold.\smallskip

In this paper, we develop an analogous deformation theory for a special class of abstract Levi-flat CR structures, the class of proper  CR submersions onto the circle (see definition \ref{properCRdef}). We deform these structure keeping the smooth submersion fixed (cf section \ref{setting}). In this setting, we prove a Kuranishi Theorem (Theorem \ref{main}), define a Kodaira-Spencer map and derive rigidity, completeness and versality criteria.\smallskip

Alternatively, our results can be seen as a study of the loop/path space of the set $\mathcal I$ of complex operators on a fixed smooth manifold $X^{diff}$. Proper CR submersions onto the circle are nothing else than deformations of the complex structures of the fiber $X^{diff}$ parametrized by the circle. So they are encoded as a loop or a path in $\mathcal I$, depending on their monodromy. From this point of view, Theorem \ref{main} describes the local action of the diffeomorphism group of the proper CR submersion onto this loop/path space, in the same way as the classical Kuranishi's Theorem describes the local action of the diffeomorphism group of a compact complex manifold onto $\mathcal I$.\smallskip

The starting point of this paper is the following easy remark. If we consider the Levi-flat CR manifold $\mathbb E_\tau\times\mathbb S^1$ (for $\mathbb E_\tau$ the elliptic curve of modulus $\tau\in\mathbb H$), then a complete space for close 
Levi-flat CR structures\footnote{Here, we consider only Levi-flat CR structures on the tangent bundle to the elliptic factor, that is we fix the differentiable type of the induced smooth foliation. Moreover, we identify two such structures if they are isomorphic through a CR isomorphism which is the identity on the $\mathbb S^1$-factor.} is given by the set of smooth maps from $\mathbb S^1$ to $\mathbb H$ close to the constant map $\tau$. It is even a local moduli space if $\tau$ is not a root of unity.\smallskip

The main result of this paper (Theorem \ref{main}) shows that this picture is still valid for Levi-flat CR manifolds which are proper submersions over the circle. Starting with any such CR manifold $\mathcal X_0$, there exists a finite-dimensional analytic space $K_{c_0}$ with a marked loop $c_0$ such that a neighborhood of $c_0$ in the loop space of $K_{c_0}$ contains all small deformations of $\mathcal X_0$. This loop space is the base of an infinite-dimensional family which is complete for $\mathcal X_0$ (Theorem \ref{completetheorem}). However, it is not always versal and we give a complete characterization of versality of our family in Theorem \ref{versaltheorem}.\smallskip

If all the fibers of the submersion $\mathcal X_0\to \mathbb S^1$ are biholomorphic, the space $K_{c_0}$ is just the common Kuranishi space of all fibers. However, if the submersion has non-isomorphic fibers, one has to build this space from different Kuranishi spaces by gluing them together. This causes a priori many technical problems. First these spaces may not have the same dimension, hence we have to fat them. Then, there is no canonical choice for the gluing maps, hence we have to make all the choices coherent. Finally, we have to make sure that the resulting glued space is really an analytic space, especially that it is Hausdorff.\smallskip

To overcome these problems, we define $K_{c_0}$ globally as a leaf space, so the problem reduces to show that this space is Hausdorff. To do that, we first have to understand Kuranishi's Theorem as giving a foliated chart of the space of complex structures. Then we show that this foliation extends to some special open sets. The leaf space of this foliation defined on some open set $\mathcal U$, when it is Hausdorff, is naturally a complete space for the structures encoded in $\mathcal U$. Our Kuranishi type moduli space is roughly speaking defined as the loop space of such a leaf space.\smallskip

It is worth emphasizing that we adopt here a completely indirect strategy, in the sense that we do not adapt the proof of the classical Kuranishi's theorem to our situation. Indeed, one of the key features of this proof is the fact that the $\bar\partial$-operator defines an elliptic resolution of the sheaf of germs of holomorphic vector fields, whose first cohomology group rules the infinitesimal deformations.  Hence,  following this direct approach amounts to finding such an elliptic resolution (e.g.  \cite{G-H-S} which proves a Kuranishi's theorem for transversely holomorphic foliations). One of its consequences is that the resulting moduli space is finite-dimensional. Now,  in the case of  a Levi flat CR structure $E$, the role of $\bar\partial$ is now taken by $\bar\partial_E$, the $\bar\partial$-operator along the leaves. The analytical properties of this operator have drawn a lot of attention, especially in connection with the exceptional minimal set problem
(see for example \cite{B-S}, \cite{H-I}, \cite{Si}). 
Nevertheless, the crucial point for deformations is that it is no longer true that it defines an elliptic complex. It is moreover easy to find examples of foliations by complex manifolds with no finite-dimensional complete family. This is even the most common case and it seems very difficult to give a useful criterion to ensure finite-dimensional completeness. For example, in the case of an irrational linear foliation of the $3$-torus by curves, finite-dimensionality is related to the arithmetic properties of the irrational slope (see \cite{Sl} or \cite{EK-S}).\smallskip

The paper is organized as follows. After giving some preliminary material in section \ref{notations}, we define and prove some useful properties of proper CR submersions over the circle in section \ref{proper}. Then we develop the general setting in section \ref{setting} and outline the results and strategy of the paper in section \ref{outline}. These two sections are crucial for the understanding of the rest of the article. The following three sections contain a thorough study of the foliated structure of the set $\mathcal I$, firstly in the neighborhood of a point (section \ref{Kuranishi}), secondly in the neighborhood of a compact set (section \ref{foliationcompact}), and thirdly in the neighborhood of a path/loop (section \ref{foliationpath}). It is the first step in the construction of the Kuranishi type moduli space and is completed in Theorem \ref{foliationtrivial}, whose technical proof is postponed to section \ref{prooftrivial}. Section \ref{Kuranishitype} contains the construction of the Kuranishi type moduli space $C^\infty (\mathbb S^1, K_{c_0})$ and the proof of Theorem \ref{main}, our main Theorem. Sections \ref{deformations}, \ref{KSmapsection}, \ref{completesection}, \ref{versalitysection}, \ref{KSbacksection} and \ref{uniquesection} develop the necessary machinery to interpret Theorem \ref{main} in terms of deformation theory in the spirit of Kodaira-Spencer and derive from it a completeness result (Theorem \ref{completetheorem}) and a versality result (Theorem \ref{versaltheorem}). Section \ref{rigiditysection} contains applications to connectedness and rigidity of the Kuranishi type moduli space. The article ends with section \ref{exemples} giving examples showing that it can be explicitely computed.

\section{Notations and preliminaries.}
\label{notations}
\noindent We begin with some notations, definitions and comments that will be used throughout the paper.
\subsection{General notations} Smooth means $C^\infty$. 
Let $X^{diff}$ be a smooth connected compact manifold of dimension $2n$. 
\vspace{5pt}\\
Given any locally trivial smooth bundle $E$ over $X^{diff}$, we denote by $\Sigma (E)$ the space of smooth sections of $E$. More generally, we denote by $\Sigma (E,B)$, or simply by $\Sigma (E)$ when the context is clear, the space of smooth sections of a locally trivial bundle $E$ with base $B$. The topology on these spaces is induced by uniform convergence on compact sets of a sequence of sections and of all its derivatives. Hence $\Sigma (E,B)$ is a Fr\'echet space. We also complete these spaces using a Sobolev norm and considering sections belonging to the Sobolev class $W_2^r$ (cf. \cite[Chapter IX]{Ku3} for more details). We denote by $\Sigma^c (E,B)$ this completion. It is a Hilbert space. We assume that $r$ is big enough to ensure that all these sections are at least of class $C^{2n+1}$. We drop any reference to this $r$ in the sequel, since it does not play any specific role.
\subsection{Classical Kodaira-Spencer theory.}
\label{classicalKSss}
The key notion is that of deformation as a flat family.
Recall the following classical definitions (cf. \cite{Su} and \cite{Ko} for additional details).

\begin{definition}
\label{deformationDef}
Let $X_0$ be a compact complex manifold. 
\begin{enumerate}
\item[(i)] An {\it analytic deformation} of $X_0$ is a flat morphism $\Pi : \mathcal X\to B$ onto a (possibly non-reduced) analytic space, together with a base-point $0\in B$ and a marking, that is a holomorphic identification $ i : X_0\to \Pi^{-1}\{0\}$.
\item[(ii)] A {\it smooth deformation} of $X_0$ is a smooth submersion $\Pi : \mathcal X\to B$ onto a smooth manifold, together with a base-point and a marking. The total space $\mathcal X$ is a endowed with a Levi-flat CR structure whose associated leaves are the level sets of $\Pi$ and the marking is assumed to be holomorphic.
\end{enumerate}
\end{definition}

In this context, Kuranishi's Theorem states the existence, for any compact complex manifold $X_0$, of an analytic deformation $\mathcal K\to K$ which satisfies
\begin{enumerate}
\item[(i)] It is {\it complete} at $0$: any analytic (resp. smooth) deformation $\mathcal X\to B$ of $X_0$ is locally isomorphic at $0$ to the pull-back of $\mathcal K$ by some analytic (resp. smooth) map $f$ from $(B,0)$ to $(K,0)$. Moreover this local isomorphism may be asked to preserve the markings.
\item[(ii)] It is {\it versal} at $0$: the (embedding) dimension of $K$ at $0$ is minimal amongst the bases of complete families for $X_0$.
\end{enumerate}
Property (ii) is known to be equivalent to the following. Given a deformation $\mathcal X\to B$, the map $f$ given by completeness is in general not unique, but its differential at $0$ is, provided only marking preserving isomorphisms are used. It can also be proven that there exists a unique germ of versal family up to isomorphism. This gives uniqueness of the Kuranishi space.\vspace{5pt}\\
From this Theorem, and from the observation that the (Zariski) tangent space at $0$ of $K$ identifies naturally with the the first cohomology group with values in the sheaf of germs of holomorphic vector fields via the so called Kodaira-Spencer map, we may derive the classical cohomological characterizations of rigidity, completeness and versality\footnote{Even if, historically, things go the other way round.}.\vspace{5pt}\\
Now, Kuranishi's Theorem has another aspect. It describes the neighborhood of a point in the set of complex operators on a fixed smooth manifold. This is indeed the point of view of Kuranishi. 

\subsection{The space of complex structures}
\label{EandIss}
We denote by $\mathcal E$, respectively $\mathcal E^c$, the set of almost-complex structures of class $C^\infty$, respectively of class $W_2^r$, on $X^{diff}$, and by $\mathcal I$, respectively $\mathcal I^c$,  the subset of $\mathcal E$, respectively of $\mathcal E^c$, formed by the integrable ones. We assume that both $\mathcal E$ and $\mathcal I$ are non-empty. Observe that our choice of $r$ implies that an element of $\mathcal I$ gives rise to a structure of complex manifold on $X^{diff}$ by Newlander-Nirenberg Theorem \cite{N-N}.
\vspace{5pt}\\
Any almost-complex operator $J$ is diagonalizable over $\mathbb C$ with eigenvalues $i$ and $-i$ and conjugated eigenspaces. This induces a splitting of the complexified tangent bundle of $X^{diff}$
\begin{equation}
\label{T}
T_{\mathbb C}X^{diff}=T\oplus \bar T
\end{equation}
Conversely, any such splitting defines a unique almost-complex structure on $X^{diff}$, which is equal (as an operator) to the multiplication by $i$ on $T$ and by $-i$ on $\bar T$. Therefore, denoting by $\text{Gr} (T_{\mathbb C} X^{diff})$ the bundle over $X^{diff}$ whose fiber at $x$ is the grassmannian of complex $n$-planes of $(T_{\mathbb C}X^{diff})_x$, we may identify $\mathcal E$ to an open subset of $\Sigma (\text{Gr} (T_{\mathbb C}X^{diff}))$. We set, more precisely
\begin{equation}
\label{E}
\mathcal E=\{T\in\Sigma (\text{Gr} (T_{\mathbb C}X^{diff}))\mid \bar T\cap T=\{0\}\}
\end{equation}
And we have
\begin{equation}
\label{I}
\mathcal I=\{T\in\mathcal E\mid [T,T]\subset T\}
\end{equation}
making $\mathcal I$ a closed subset of $\mathcal E$.
\vspace{5pt}\\
Since $\text{Gr} (T_{\mathbb C}X^{diff})$ is a complex bundle, its space of $W_2^r$-sections is a Hilbert manifold over $\mathbb C$ with tangent space at some section $\sigma$ equal to the Hilbert space of $W_2^r$-sections of the complex vector bundle $ \sigma^*T\text{Gr}(T_{\mathbb C}X^{diff})$. In particular, 
the open subset $\mathcal E^c$ is a $\mathbb C$-Hilbert manifold. In the same way, $\mathcal E$ is a Fr\'echet manifold \cite{Hamilton}. \vspace{5pt}\\
We will make use of the formalism of Hilbert analytic spaces as defined in \cite{Do1} and \cite{Do2}. 
The subset $\mathcal I^c$ is a $\mathbb C$-analytic Hilbert space (cf. \cite{Do1}). By abuse of notations, we will say that $\mathcal I$ is a {\it Fr\'echet analytic space}. More generally, we say that a closed subset $I$ of a Fr\'echet manifold modeled on some Fr\'echet space $\Sigma (E,B)$ is a {\it Fr\'echet analytic space} if it is locally the restriction of some Hilbert analytic set of  $\Sigma^c (E,B)$.\vspace{5pt}\\
Given $J\in\mathcal E$ (respectively $J\in\mathcal I$), we denote by $X_J$ the almost-complex (respectively complex) manifold $(X^{diff}, J)$. We also set $X_0$ for $X_{J_0}$.\vspace{5pt}\\
Deformations can also be defined as analytic (resp. smooth) families of complex operators. To be more precise,

\begin{definition}
Let $B$ be an analytic space (resp. a smooth manifold) with base-point $0$. Consider a family $(J_t)_{t\in B}$ of elements of $\mathcal I$. We say that $(J_t)_{t\in B}$ is {\it analytic} (resp. {\it smooth}) if endowing each fiber $X^{diff}\times \{t\}$ of the projection $X^{diff}\times B\to B$ with the complex structure $J_t$ turns it into an analytic (resp. smooth) deformation of $X$.
\end{definition}

The Kuranishi space $K$ of $X$ defines such a family of complex operators. It follows from \cite[Theorem 8.1]{Ku3} that this family is analytic. Hence the Kuranishi space $K$ of $X$ naturally defines an analytic deformation $\Pi : \mathcal K\to K$ once chosen a marking. It is called the {\it Kuranishi family}. From the existence of the Kuranishi family, we obtain the following. Given an analytic (resp. smooth) map $f : B\to \mathcal I$, the family $(f(t))_{t\in B}$ is analytic (resp. smooth). 

\subsection{The diffeomorphism group and its action}
\label{Diffss}
Let $ \text{Diff} (X^{diff})$, respectively $\text{Diff}^c (X^{diff})$, be the group of diffeomorphisms of class $C^\infty$, respectively of class $W_2^{r+1}$, of $X^{diff}$ and let $ \text{Diff}^0 (X^{diff})$, respectively $\text{Diff}^{0,c} (X^{diff})$  be the subgroup of elements smoothly isotopic to the identity. They act on $\mathcal E$ as follows: for $J$ an almost complex operator and $f\in \text{Diff}(X^{diff})$, we have
\begin{equation}
\label{action}
(f\cdot J)_x (v)=d_{f^{-1}(x)}f\circ J_{f^{-1}(x)}\circ (d_xf)^{-1} (v)
\end{equation}
for $(x,v)\in TX^{diff}$, or, for $T\in\mathcal E$ in the presentation \eqref{E},
$$
(f\cdot T)_x=(d_{f^{-1}(x)}f)(T_{f^{-1}(x)}).
$$
Besides, by definition, both preserve the almost-complex structures, that is $f$ realizes an isomorphism between $X_J$ and $X_{f\cdot J}$.\vspace{5pt}\\
The group $\text{Diff}^c (X^{diff})$ being an open subset of the Hilbert manifold of $W_2^{r+1}$ maps from $X^{diff}$ to the complex manifold $X_J$, is also a Hilbert manifold and the action is $\mathbb C$-analytic locally at $J$ (cf. \cite{Do1}). The group $ \text{Diff} (X^{diff})$ is a Fr\'echet Lie group.\vspace{5pt}\\
Let $W$ be a neighborhood of the identity in $\text{Diff} (X^{diff})$. Observe that, if $W$ is small enough, every $\phi\in W$ can be constructed as follows. Endow $X^{diff}$ with a smooth riemannian metric. There exists a smooth vector field $\xi$ close to $0$ such that the map
\begin{equation}
\label{exp}
x\in X^{diff}\buildrel e(\xi) \over \longmapsto \gamma_{x,\xi(x)}(1)
\end{equation}
is exactly $\phi$. Here 
$$
\gamma_{x,\xi(x)}\ :\ \mathbb R^+\longrightarrow X^{diff}
$$
is the geodesic starting at $x$ with initial velocity $\xi(x)$.
\vspace{5pt}\\
Conversely, there exists an open neighborhood V of $0$ in $\Sigma (TX^{diff})$ such that, for every $\xi\in V$, the map $e(\xi)$ defined by \eqref{exp} is a diffeomorphism of $X^{diff}$. Hence we constructed in that way a chart $e$ from $V$ to $W$.
\vspace{5pt}\\
We denote by $\text{Aut}(X_J)$, respectively $\text{Aut}^0(X_J)$, the finite-dimensional Lie group of biholomorphic transformations of the complex manifold $X_J$, respectively the connected component of the identity in this group. Observe that $\text{Aut}(X_J)$ is the isotropy group of the action \eqref{action} at $J$.
\subsection{Almost complex preserving maps}
\label{acss}
Given two subsets $U$ and $V$ of $\mathcal E$ with $U$ open and a smooth map $F$ from $U$ to $V$, we say that $F$ is {\it almost-complex preserving} if, for each $J\in U$, the manifolds $X_J$ and $X_{F(J)}$ are CR isomorphic. We denote
$$
\begin{CD}
U @ > F> a.c. > V
\end{CD}
$$
Notice that an a.c. map is a special type of equivariant map. Indeed, the set of a.c. maps corresponds exactly to the set of equivariant maps which descends as the identity on the quotient space $\mathcal E/\text{Diff}(X^{diff})$. We extend this notion to the following cases. Firstly, letting $U$ and $V$ as before and letting $W$ be an open subset of a topological $\mathbb C$-vector space, we say that a smooth map $F$ from $U$ to $V\times W$ (respectively from $U\times W$ to $V$) is {\it almost-complex preserving} if the composition
$$
\begin{CD}
U @ >F>> V\times W@ >\text{1st projection} >> V
\end{CD}
$$
respectively
$$
\begin{CD}
U@ >\text{inclusion}>>U\times \{w\}\subset U\times W@ >>> V
\end{CD}
$$
is almost-complex preserving (respectively almost-complex preserving for all $w\in W$). Secondly, if $F$ maps smoothly some open set $U$ of $\mathcal E$ to some manifold or analytic space $K$, then we say that $F$ is {\it almost-complex preserving} if we have
$$
F(J)=F(J')\quad\Longrightarrow\quad X_J\text{,}\ X_{J'} \text{ and }X_{F(J)}\text{ are CR isomorphic}.
$$
Observe that every point of $K$ is thus associated to an isomorphism class of complex structures.
\subsection{Paths and loops}
All paths and loops are smooth. Given $M$ a finite dimensional or a Fr\'echet manifold, respectively analytic space, the loop space $C^\infty(\mathbb S^1, M)$ is topologized by uniform convergence of sequences of maps $\mathbb S^1\to M$ and of all their derivatives. It is a Fr\'echet manifold, respectively analytic space. We say that two such loops are {\it close} if they belong to a small open set of  $C^\infty(\mathbb S^1, M)$. Similar definition holds for paths.

\section{Proper CR submersions.}
\label{proper}

Let us define the objects we want to deform and study.

\begin{definition}
\label{properCRdef}
A proper smooth submersion $\pi : \mathcal X\to \mathbb S^1$ is called a {\it a proper CR submersion} if $\mathcal X$ is endowed with a Levi-flat integrable CR structure which is tangent to the fibers of $\pi$.
\end{definition}

As a smooth manifold, a proper CR submersion is a locally trivial smooth fiber bundle over the circle, thanks to Ehresmann's Lemma. The fiber, that we denote by $X^{diff}$, is a smooth compact manifold. We assume that it is connected. In other words, $\mathcal X$ is diffeomorphic to 
\begin{equation}
\label{CRdiffmodel}
X_\phi :=(X^{diff}\times [0,1])/\sim
\quad\text{ where }\quad (x,0)\sim (x',1)\iff x'=\phi(x).
\end{equation}
Here $\phi$ is a fixed diffeomorphism of $X^{diff}$, classically called the {\it monodromy} of $X_\phi$. Recall also that $X_\phi$ and $X_{\phi'}$ are diffeomorphic if $\phi$ and $\phi'$ are isotopic.\smallskip

As a CR manifold, each fiber of $\mathcal X$ is a copy of $X^{diff}$ equipped with a complex structure. The only difference between a proper CR submersion and a smooth deformation over the circle (cf. definition \ref{deformationDef}) is that here there is no marked point.
\vspace{5pt}\\
By Fischer-Grauert's Theorem (see \cite{Me2} for the version we use), if all the fibers of a proper CR submersion are biholomorphic to a fixed manifold $X_0$, then it  is locally trivial, that is it satisfies
\begin{equation}
\label{CRtrivial}
\begin{CD}
\pi^{-1}(U)@ >\text{CR isomorphism}>>U\times X_0\cr
@ V\pi VV @ VV\text{1st projection}V\cr
U@ >>\text{Identity}> U
\end{CD}
\end{equation}
in a neighborhood $U$ of any point of the circle. 

\begin{definition}
\label{CRbundle}
We call a proper CR submersion a {\it CR bundle} if it satisfies \eqref{CRtrivial}. It is {\it trivial} if it is globally CR isomorphic to $X_0\times\mathbb S^1$;  {\it non-trivial} otherwise. 
\end{definition}

We have immediately

\begin{lemma}
\label{CRtriviallemma}
A CR bundle $\pi : \mathcal X\to\mathbb S^1$ is trivial if and only if its monodromy $\phi$ is isotopic to the identity.
\end{lemma}

Of course, in the general case,  the fibers of a proper CR submersion have distinct complex structures. Indeed, 
choosing a smooth model \eqref{CRdiffmodel} for $\mathcal X$, we may identify it as a CR manifold with a smooth path $c$ in $\mathcal I$ such that
\begin{equation}
\label{cends}
c(1)=\phi\cdot c(0)
\end{equation}
and
\begin{equation}
\label{cpath}
c_\phi : t\in (-\epsilon, \epsilon)\longmapsto \left \{
\begin{aligned}
&c(1+t) \quad\text{ if } t\leq 0\cr
&\phi\cdot c(t)\quad\text{ if }t\geq 0
\end{aligned}\right .
\qquad\text{ is smooth.}
\end{equation}
that is, $\mathcal X$ is CR isomorphic to the smooth submersion
\begin{equation}
\label{smoothCR}
\pi^{diff} : \mathcal X^{diff}\to\mathbb S^1
\end{equation}
endowed with the family of complex operators 
\begin{equation}
\label{CRfamily}
(\pi^{-1}(\exp {2i\pi t}), c(t))_{\exp {2i\pi t}\in \mathbb S^1},
\end{equation}
the identification of the endpoints being realized thanks to \eqref{cends} and \eqref{cpath}.
\vspace{5pt}\\
Hence, we may and will write $\mathcal X_0=(\mathcal X^{diff}, c_0)$ and $\mathcal X_c=(\mathcal X^{diff}, c)$ exactly in the same way as we write $X_0=(X,J_0)$ and $X_J=(X,J)$ for compact complex manifolds.
\vspace{5pt}\\
Of course, as in the case of complex structures, the encoding of a proper CR submersion as a loop $c$ is far from being unique, since we can move $c$ in $\mathcal I$ using the $\text{Diff}^0(X^{diff})$ action on each point $c(t)$.

\section{Setting.}
\label{setting}

In this section, we develop the general setting of the paper. 
\vspace{5pt}\\
Our aim is to study the deformations of proper CR submersions with fixed differentiable type of the induced smooth foliation. Precise definitions of such deformations in the spirit of Kodaira-Spencer are given in section \ref{deformations}. For the moment, let us just say that we keep the smooth model \eqref{CRdiffmodel} and we modify the family \eqref{CRfamily} of complex operators along the fibers of the submersion onto the circle. Moreover, we identify two such structures if they are isomorphic through a CR isomorphism which is the identity on the $\mathbb S^1$-factor. 
\vspace{5pt}\\
In other words, we encode a proper CR submersion as a path $c_0$  in $\mathcal I$ and a monodromy map $\phi$ satisfying \eqref{cends}, \eqref{cpath}.  Any small deformation of it is given by a path close to $c_0$ up to the action of the diffeomorphisms.
\vspace{5pt}\\
Given a compact complex manifold $X_0=(X^{diff}, J_0)$, the classical Kuranishi's Theorem constructs an analytic space that describes all small deformations of the complex structure of $X_0$ in a minimal way. Its proof consists in giving a precise description of a neighborhood of $J_0$ in $\mathcal I$ and of the action \eqref{action} onto it.
\vspace{5pt}\\
In the same way, given a proper CR submersion $\mathcal X=(\mathcal X^{diff}, c_0)$, in order to construct an analytic space that describes all its small deformations in a minimal way, we need to understand the neighborhood of the path $c_0$ in $\mathcal I$ and the action of the diffeomorphisms onto it.
\vspace{5pt}\\
In the next two subsections, we describe the neighborhood of the path $c_0$ in $\mathcal I$ as the loop space of an infinite-dimensional Fr\'echet analytic space; and the action of the diffeomorphisms as the action of a well defined Fr\'echet Lie group. These subsections are built in complete analogy with subsections \ref{EandIss} and \ref{Diffss} to which we refer.

\subsection{The set of proper CR submersions}
\label{CRsubss}
As proved in the section \ref{proper}, the space of CR submersions compatible with $\pi^{diff}$ is just the set
\begin{equation}
\label{IXdiff}
\mathcal I(\mathcal X^{diff})=\{c : [0,1]\to \mathcal I\mid c\text{ is smooth and satisfies }\eqref{cends}, \eqref{cpath}\}
\end{equation}
and we may define similarly the space of {\it almost CR submersions} as the set
\begin{equation}
\label{EXdiff}
\mathcal E(\mathcal X^{diff})=\{c : [0,1]\to \mathcal E\mid c\text{ is smooth and satisfies } \eqref{cends}, \eqref{cpath}\}.
\end{equation}
If the monodromy of $\mathcal X_0$ is the identity, then the path $c_0$ is a loop and we have
\begin{lemma}
\label{loopstructureIlemma}
We have
\begin{enumerate}
\item[(i)] The space $\mathcal E(\mathcal X^{diff})$ is the loop space $C^{\infty}(\mathbb S^1,\mathcal E)$, hence is a Fr\'echet manifold.
\item[(ii)] The space $\mathcal I(\mathcal X^{diff})$ is the loop space $C^{\infty}(\mathbb S^1,\mathcal I)$, hence is a Fr\'echet analytic space.
\end{enumerate}
\end{lemma}

In particular, if $U_\phi$ is a neighborhood of $c_0$ in $\mathcal E$, then $C^{\infty}(\mathbb S^1,\mathcal I\cap U_\phi)$ contains all small deformations of $\mathcal X_0$.
\vspace{5pt}\\
If the monodromy of $\mathcal X_0$ is not the identity, then we have to twist the construction, gluing elements $J\in\mathcal E$ to $\phi\cdot J$. Of course this cannot be done globally on $\mathcal E$, since we want the resulting quotient to still being a Fr\'echet manifold. But it can be done locally, in a neighborhood of $c_0$.
\vspace{5pt}\\
We take a connected neighborhood $U$ of (the image of) $c_0$ and decompose it into a union of connected open subsets
\begin{equation}
\label{adaptedI}
U=U_1\cup\hdots\cup U_k
\end{equation}

\begin{definition}
\label{adapteddef}
We say that the covering {\rm \eqref{adaptedI}} is {\it adapted} if it satisfies the following conditions.
\begin{enumerate}
\item[(i)] Each open set $U_i$ is the domain of a Fr\'echet chart.
\item[(ii)] If $c_0$ is a path with distinct endpoints, then an intersection $\overline{U_i}\cap \overline{U_j}$ is non empty if and only if $j=i$ or $j=i+1$ or $i=j+1$. In this case, it is connected.
\item[(iii)] If $c_0$ is a loop, then an intersection $\overline{U_i}\cap \overline{U_j}$ is non empty if and only if $j=i$ or $j=i+1$ or $i=j+1$ or $(i,j)=(1,k)$ or $(j,i)=(1,k)$. In this case, it is connected.
\item[(iv)] The path $c_0$ intersects each $U_i$ and each non empty intersection $U_i\cap U_j$ in a connected path.
\end{enumerate}
\end{definition}

Observe that $c_0$ admits an adapted covering as soon as it has no self-intersection. Assume that $\phi$ is not an automorphism of $X_{c_0(0)}$. The analogous statement to Lemma \ref{loopstructureIlemma} is

\begin{lemma} 
\label{loopstructureIIlemma}
Assume that $c_0$ does not self-intersect and that $c'_0(0)$ is different from $0$. For any sufficiently small neighborhood $V$ of $c_0$ in $\mathcal E(\mathcal X^{diff})$, there exists a Fr\'echet manifold $U_\phi$ and a Fr\'echet analytic space, that we denote by abuse of notation $U_\phi\cap\mathcal I$, such that
\begin{enumerate}
\item[(i)] The loop space $C^\infty(\mathbb S^1, U_\phi)$ is a Fr\'echet manifold homeomorphic to $V$.
\item[(ii)] The loop space $C^\infty(\mathbb S^1, U_\phi\cap \mathcal I)$ is a Fr\'echet analytic space homeomorphic to $V\cap\mathcal I(\mathcal X^{diff})$.
\end{enumerate}
\end{lemma}

\begin{proof}
This is a typical example of the problems we will meet in the proof of say, Theorem \ref{main}. So even if it is a much easier case, we give a detailed treatment. \vspace{5pt}\\
Choose an adapted covering $U$ of $c_0$ with $k>2$. Our assumption that $\phi$ is not a biholomorphism of $X_{c_0(0)}$ implies that $c_0(1)$ is different from $c_0(0)$. We assume that $\phi$ sends bihomorphically $U_1$ onto $U_k$.  We define now $U_\phi$ as the set obtained from $U$ by identifying each point $J$ of $U_1$ with the corresponding point $\phi\cdot J$ of $U_k$. We put on $U_\phi$ the quotient topology. 
\vspace{5pt}\\
Observe that we may assume it is Hausdorff, shrinking $U$ if necessary. To check that, we switch to the Hilbert setting. To simplify notations, we still denote by $U$ and $U_i$ the corresponding open subsets of $\mathcal E^c$. Now, the only problem that could appear is the following. If we can find a sequence $(x_n)$ in $U_1$
such that
\begin{equation}
\label{NHx}
\lim x_n\in U\cap \overline{U_1}\setminus U_1
\end{equation}
and 
\begin{equation}
\label{NHy}
y:=\lim \phi (x_n)\in U\cap \overline{U_k}\setminus U_k
\end{equation}
then $U_\phi$ is not separated at $y$. This can be avoided as follows.
\begin{figure}[H]
    \centering
    \includegraphics[width=.6\textwidth]{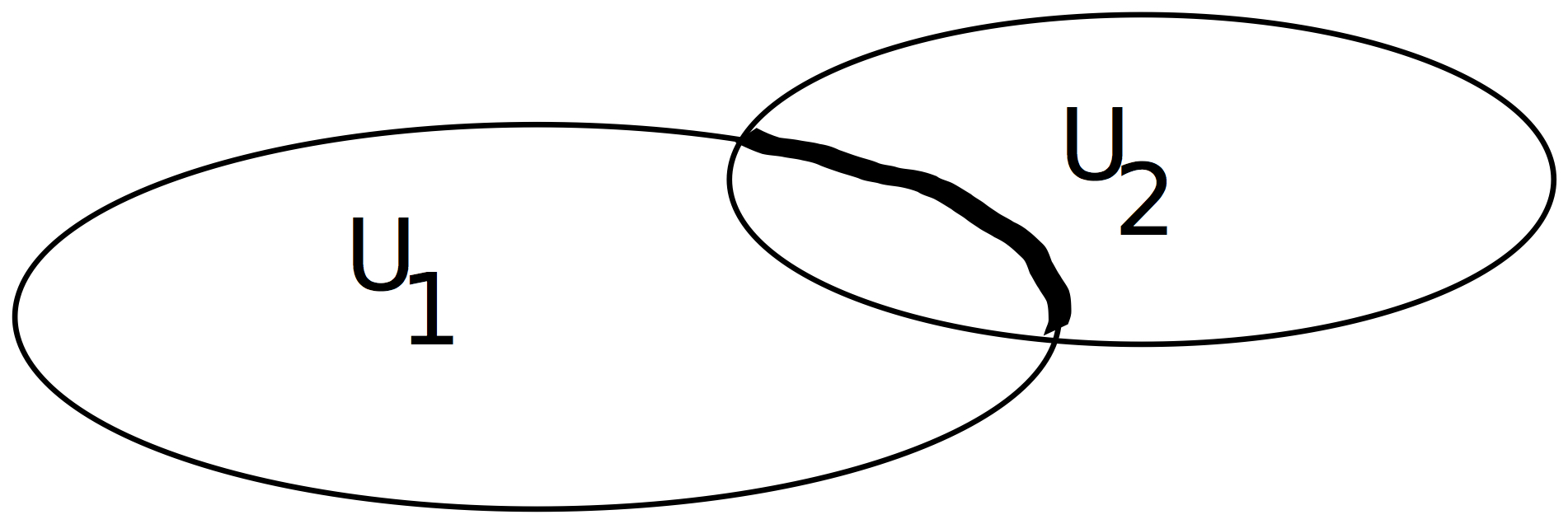}
    \caption{$U\cap \overline{U_1}\setminus U_1$.}
    \label{U1}
\end{figure}
As $c_0$ maps into a Hilbert manifold and $c'_0(0)\not =0$, we may define
\begin{equation}
\label{Hperp}
H:=(c'_0(0))^\perp
\end{equation}
We then extend $c_0$ to $(-\epsilon, 0]$ for $\epsilon$ small enough so that on $(-\epsilon,\epsilon)$
$$
\begin{aligned}
\langle \overrightarrow{c_0(0)c_0(t)}, c_0'(0)\rangle &>0\iff 0<t<\epsilon\cr
&<0\iff -\epsilon<t<0
\end{aligned}
$$
Shrinking $U$ if necessary, we may assume that $U_1$ is a ball and that $H$ cuts $U_1$ into two connected components, say $H^+$ and $H^-$. Moreover, we assume that $U_2\subset H^+$ and, using $U_k=\phi\cdot U_1$ and \eqref{cpath} and $c'_0(0)\not =0$, we have that $U_{k-1}\subset \phi\cdot H^-$.\vspace{5pt}\\
Let now $x$ be a point of $U\cap \overline{U_1}\setminus U_1$. Then $x$ belongs to $U_2$, hence to $H^+$. As a consequence, if $(x_n)$ converges to $x$, then  $\phi(x_n)$ belongs to $\phi\cdot H^+$ for $n$ big enough. Since $U\cap \overline{U_k}\setminus U_k$ is included in $\phi\cdot H^-$, we see that it is not possible to find a sequence $(x_n)$ satisfying \eqref{NHx} and \eqref{NHy}.
\vspace{5pt}\\
Hence $U_\phi$ is Hausdorff. Since it is obtained from the Hilbert manifold $U$ by an open gluing, it is also a Hilbert manifold. We may thus speak of smooth maps into $U_\phi$.
\vspace{5pt}\\
If $\mathcal X$ is a CR submersion close to $\mathcal X_0$, it is represented by a smooth path $c$ in $U$ whose endpoints belong to $U_1$ and $U_k$ respectively, and which satisfies \eqref{cends} and \eqref{cpath}. Hence it is represented by a smooth loop in $U_\phi$. Reciprocally it is clear that any loop in $U_\phi$ formed by integrable structures lifts to a point of $\mathcal I(\mathcal X^{diff})$ close to $c_0$, thus defining a CR submersion close to $\mathcal X_0$.
\vspace{5pt}\\
In other words, denoting abusively $U_\phi\cap\mathcal I$ the subset of points of $U_\phi$ corresponding to integrable structures, we may take the space $C^\infty (\mathbb S^1, U_\phi)$ as a neighborhood of $c_0$ in $\mathcal E(\mathcal X^{diff})$; and the space $C^\infty (\mathbb S^1, U_\phi \cap\mathcal I)$ as a neighborhood of $c_0$ in $\mathcal I(\mathcal X^{diff})$.
\end{proof}

If $\phi$ is a non-trivial automorphism of $X_{c_0(0)}$, then the previous construction does not work. For example, if $c_0$ is a constant path, then it would consist in taking the quotient of a neighborhood of $c_0(0)$ in $\mathcal I$ by the action generated by $\phi$, which fixes $c_0(0)$. In the sequel, we will avoid this type of paths and monodromy by using the action of the diffeomorphisms. This will also allow us to fulfill the hypotheses of Lemma \ref{loopstructureIIlemma}.

\subsection{The diffeomorphism group and its action}
\label{CRDiffss}
Let $\text{Diff}_\pi (\mathcal X^{diff})$ be the group of bundle isomorphisms of $\pi^{diff}$ which descend as the identity on the base $\mathbb S^1$.  We will focus indeed to its subgroup $\text{Diff}^0_\pi (\mathcal X^{diff})$ of elements isotopic to the identity in $\text{Diff}_\pi (\mathcal X^{diff})$. We can describe more precisely this last group. \vspace{5pt}\\
As usual, take a presentation \eqref{CRdiffmodel}. Then $\text{\rm Diff}^0_\pi (\mathcal X^{diff})$ is homeomorphic to the set of paths
\begin{equation}
\label{diffpistructure}
\{F : [0,1]\to \text{\rm Diff}^0(X^{diff})\mid \phi\circ F(0)=F(1)\circ\phi\}
\end{equation}
satisfying moreover 
\begin{equation}
\label{diffpath}
F_\phi : t\in (-\epsilon, \epsilon)\longmapsto \left \{
\begin{aligned}
&F(1+t)\circ \phi \quad\text{ if } t\leq 0\cr
&\phi\circ F(t)\quad\text{ if }t\geq 0
\end{aligned}\right .
\qquad\text{ is smooth.}
\end{equation}
We may think of its elements as an isotopy
$$
\left (f_t\ :\ X^{diff}\longrightarrow X^{diff}\right )_{e^{2i\pi t}\in\mathbb S^1}
$$
the identification at the endpoints using \eqref{diffpistructure} and \eqref{diffpath}. If $\phi$ is the identity, then $\text{Diff}^0_\pi (\mathcal X^{diff})$ is the loop space $C^\infty(\mathbb S^1,\text{Diff}^0(X^{diff}))$. This remark should be compared with Lemma \ref{loopstructureIlemma}. In the same way, we could prove a result similar as Lemma \ref{loopstructureIIlemma} in the case $\phi$ is not an automorphism of $X_{c_0(0)}$. However, we will not need such a statement, but rather will use the following
\begin{lemma}
\label{Diffstructurelemma}
The group $\text{\rm Diff}^0_\pi (\mathcal X^{diff})$ is a Fr\'echet Lie group.
\end{lemma}

\begin{proof}
It is enough to construct a Fr\'echet chart at the identity. The fact that the group operations are smooth comes from the fact that $\text{Diff}^0(\mathcal X^{diff})$ is a Fr\'echet Lie group cf. \cite{Hamilton}.\vspace{5pt}\\
Consider the set $\Sigma_\pi (T\mathcal X^{diff})$ of smooth vector fields of $\mathcal X^{diff}$ tangent to the fibers of $\pi$. We have
$$
\Sigma_\pi (T\mathcal X^{diff})\simeq \{\sigma : t\in\mathbb S^1\longmapsto \sigma_t\in \Sigma(T\pi^{-1}(t))\}.
$$
Taking a riemannian metric along the fibers, we may define, as in \eqref{exp}, a map
\begin{equation}
\label{exppi}
e_\pi : \Sigma_\pi (T\mathcal X^{diff})\longrightarrow \text{Diff}_\pi (\mathcal X^{diff})
\end{equation}
which satisfies, for all $t\in\mathbb S^1$,
$$
(e_\pi)_t : \xi\in \Sigma (TX^{diff})\simeq\Sigma (T\pi^{-1}\{t\})\longmapsto e(\xi)\in \text{Diff}(X^{diff})\simeq \text{Diff}(\pi^{-1}\{t\})
$$
As in the pointwise case, $e_\pi$ realizes a diffeomorphism between an open neighborhood of the zero-section in $\Sigma_\pi (T\mathcal X^{diff})$ and a neighborhood of the identity in $\text{Diff}_\pi (\mathcal X^{diff})$.
\end{proof}

The group $\text{Diff}_\pi (\mathcal X^{diff})$ acts on $C^\infty (\mathbb S^1, U_\phi)$ in the obvious way. That is, given $c\in C^\infty (\mathbb S^1, U_\phi)$ and $f\in \text{Diff}_\pi (\mathcal X^{diff})$, we have
\begin{equation}
\label{piaction}
(f\cdot c)_t=f_t\cdot c_t
\qquad \text{ for all }t\in\mathbb S^1
\end{equation}
where the action on the right-hand side is defined in \eqref{action}. Indeed, \eqref{piaction} plays the role of \eqref{action} for proper CR submersions.\vspace{5pt}\\
 In other words, 
 \begin{lemma}
 \label{pathequivalencelemma}
 The following two statements are equivalent
 \begin{enumerate}
 \item[(i)] Two paths $c\in C^\infty (\mathbb S^1, U_\phi)$ and $c'\in C^\infty (\mathbb S^1, U_\phi)$ encode the same proper CR submersion up to CR isomorphisms inducing the identity on the base $\mathbb S^1$.
 \item[(ii)] We have $c'=f\cdot c$ for some $f$ in $\text{\rm Diff}_\pi (\mathcal X^{diff})$.
 \end{enumerate}
 \end{lemma}
Before finishing this subsection, let us  introduce some hypotheses about the paths we use. We will always assume without loss of generality that
\begin{hypothesis}
\label{noselfint}
We have
\begin{enumerate}
\item[(i)] The path $c_0$ has no self-intersection in $\mathcal I$.
\item[(ii)] If the monodromy of $\mathcal X_0$ is isotopic to the identity, we take $c_0$ to be a loop and $\phi$ to be the identity. 
\item[(iii)] If $\mathcal X_0$ is a trivial CR bundle, then we take $c_0$ to be a constant loop.
\item[(iv)] If $\mathcal X_0$ is a non-trivial CR bundle, then we take $c_0$ to be a non constant path with distinct endpoints with $c'_0(0)\not =0$.
\end{enumerate}
\end{hypothesis}

As for (i), in case $c_0$ does self-intersect, just move it locally along the $\text{Diff}^0_\pi (\mathcal X^{diff})$-orbits in a finite number of points to destroy the self-inter\-sec\-tions. As for (iv), use also the action so that $c_0$ is a path with distinct endpoints in a single $\text{Diff}^0_\pi (\mathcal X^{diff})$-orbit. The interest of (iv) should be clear when considering Lemma \ref{loopstructureIIlemma} and the discussion below.
\subsection{Almost-complex preserving maps}
\label{acdiffmapsss}
Given two subsets $U$ and $V$ of $\mathcal E (\mathcal X^{diff})$ with $U$ open, $W$ an open subset of a topological $\mathbb C$-vector space, and an analytic map $F$ from $U$ to $V\times C^\infty (\mathbb S^1,W)$, we say that $F$ is {\it almost-complex preserving}  if the composition
$$
\begin{CD}
c\in U @ >F>> V\times C^\infty (\mathbb S^1,W)@ >\text{1st projection} >> G(c)\in V
\end{CD}
$$
is almost-complex preserving for each $t\in\mathbb S^1$, that is the complex manifolds $X_{c(t)}$ and $X_{G(c(t))}$ are isomorphic for each $t$. 
\vspace{5pt}\\
Thus we extend the notion of a.c. maps to $\mathcal E(\mathcal X^{diff})$. Observe that an a.c. map is equivariant with respect to the action of $\text{Diff}^0_\pi (\mathcal X^{diff})$.

\section{Outline of the paper.}
\label{outline}
In this section, we outline the results and the strategy of the paper. The main Theorem is Theorem \ref{main}. It is a generalization of the classical Kuranishi's Theorem to proper CR submersions. From this result, we derive rigidity, completeness and versality criteria. Theorem \ref{main} states that the set of proper CR submersions close to a fixed one, described in section \ref{setting} as the loop space of the infinite-dimensional Fr\'echet analytic space $U_\phi\cap\mathcal I$, can also be described as the loop space of a finite-dimensional analytic space $K_{c_0}$. The central idea of the proof is to show firstly that the neighborhood $U_\phi\cap\mathcal I$ is foliated with leaves included in the orbits of $\text{Diff}^0(X^{diff})$; and secondly that this foliation has a Hausdorff leaf space $K_{c_0}$. Hence $U_\phi\cap\mathcal I$ admits an a.c. retraction onto $K_{c_0}$. As a consequence every loop in $U_\phi\cap\mathcal I$ can be pushed to a loop of $K_{c_0}$ and small deformations are also described by the loop space of the finite-dimensional $K_{c_0}$. \smallskip

Let us be more specific. For that, it is important to go back once again to the classical Kuranishi's Theorem. It can be restated as: there exists a retraction along the orbits of $\text{Diff}^0(X^{diff})$ of a neighborhood of $J_0$ in $\mathcal I$ onto a finite-dimensional analytic subspace called the Kuranishi space of $J_0$. Actually it is shown that there exists a local holomorphic foliation of $\mathcal I$ at $J_0$ whose leaves are modeled on $\text{Diff}^0(X^{diff})/\text{Aut}^0(X_0)$. The Kuranishi space of $X_0$ embeds as a transverse section and the neighborhood of $J_0$ retracts onto this transverse section following the leaves. This is what we explain in section \ref{Kuranishi}.\smallskip

The above mentioned foliation of $U_\phi\cap\mathcal I$ will be obtained by extending this local foliated chart into a foliated atlas of $U_\phi\cap\mathcal I$. Indeed, take any compact set $C$ containing $J_0$. Then Theorem \ref{foliation} states that the local foliation at $J_0$ extends as a foliation defined in a neighborhood of $C$. Once again, the leaves are modeled on $\text{Diff}^0(X^{diff})/\text{Aut}^0(X_0)$ and the Kuranishi space of $X_0$ can be realized as a transverse section at $J_0$. But at another point $J$, the Kuranishi space of $X_J$ does not always identify with a transverse section at $J$. Of course, this annoying point comes from the fact that the leaves are modeled on $\text{Diff}^0(X^{diff})/\text{Aut}^0(X_0)$, and since the automorphism group of $X_J$ depends on $J$, this space may be different from $\text{Diff}^0(X^{diff})/\text{Aut}^0(X_J)$. We call this foliation the $L$-foliation, where $L$ is a fixed complementary subspace in the space of vector fields to the Lie algebra of $\text{Aut}^0(X_J)$.\smallskip

We come to the most important property of this foliation, stated as Theorem \ref{foliationtrivial}, the main result of this paper. Take $C$ to be the path $c_0$ - and $J_0$ a point of $C$. Then, the leaf space of the $L$-foliation of a neighborhood of $C$ is a Hausdorff analytic space $K_{c_0}$. The line of arguments of the proof is very classical in (finite-dimensional) foliation theory. We assume Hypothesis \eqref{noselfint} and we take an adapted covering \eqref{adaptedI} constituted by foliated charts. By induction on $k$, the number of charts in the adapted covering, we prove that the local transverse sections can be locally pushed along the leaves to match and form a global transverse section to the foliation. Of course, such local movements use a partition of unity, so this global transverse is only a smooth object. In other words, the foliation is smoothly trivial, not holomorphically trivial. However, this Hausdorff section can be turned into an analytic space by making use of the foliated atlas: it is locally modeled onto the analytic local section and changes of charts are given by the transverse component of the changes of charts of the foliated atlas hence are analytic. Observe that this classic trick appears in a completely different setting as a method for constructing complex non K\"ahler manifolds, see \cite{Mathese} and the famous Bogomolov conjecture \cite{Bog}, \cite{DG}. In any case, since we are in an infinite-dimensional case, and with singular transverse sections, we have to check carefully that the classical proof can be adapted to this situation.  The resulting demonstration is quite technical and is postponed in section \ref{prooftrivial}. \smallskip

As a consequence of Theorem \ref{foliationtrivial}, every proper CR submersion close to $\pi : \mathcal X\to \mathbb S^1$ can be encoded as a path into $K_{c_0}$. Hence we end with the loop space of $K_{c_0}$ playing the role of the Kuranishi space in the classical setting. We call this loop space a Kuranishi type moduli space and denote it by $K^g$. The retraction along the $\text{Diff}^0$-orbits from a neighborhood of $J_0$ onto its Kuranishi space is replaced with the retraction along the $\text{Diff}^0$-orbits of the loop space of $U_\phi$ onto $K^g$. \smallskip

Next, we explore the deformation theory of proper CR submersions in the spirit of Kodaira-Spencer. This supposes to define a notion of deformation of proper CR submersions as flat family. This is what we do in section \ref{deformations}, before defining a Kodaira-Spencer map in section \ref{KSmapsection}. For this, we follow the original approach by Kodaira-Spencer, writing down a short exact sequence of sheaves which encodes the obstructions for lifting vector fields of the base space as infinitesimal automorphisms of the family. We then can construct a deformation $\mathcal K^g\to K^g$ of $c_0$ from Theorem \ref{main} and prove that this family is complete in Theorem \ref{completetheorem}. This is quite a straightforward result since we develop before the necessary definitions and tools. The case of versality is much involved. We give a criterion for versality and another one for universality in Theorem \ref{versaltheorem}. The idea behind these criteria is simple. We just explain the criterion for versality. Basically the analytic space $K_{c_0}$ is constructed from the Kuranishi spaces of the fibers of the CR submersion $\mathcal X_0$. Assume firstly that the monodromy of $\mathcal X_0$ is isotopic to the identity. If the dimension of these Kuranishi spaces does not change along the circle, then $K_{c_0}$ is roughly speaking the union of a finite number of these spaces glued along their common intersections and the loop space of $K_{c_0}$ is versal. However, if the dimension varies, $K_{c_0}$ must have the maximal dimension and versality is lost. Notice also that we speak of the dimension of the regular part. This explains why criterion \eqref{versalcondition} is in terms of the dimension of the automorphism group and not of the first cohomology group with values in the tangent sheaf. Secondly, in the case of a non-trivial CR bundle, our construction always yields a non versal Kuranishi type moduli space. This comes from the fact that the natural candidate for $K_{c_0}$ would be the quotient of the Kuranishi space of the single fiber (since we are in the CR bundle case, there is a single fiber up to biholomorphism) by the automorphism serving as monodromy. However, in many cases this is not an analytic space but an analytic stack. But, more deeply, even in the case where such a quotient is an analytic space, the loop space of this quotient would not be complete, since the non-trivial CR bundle should correspond via pull-back to a constant loop, see example \ref{badcaseexample}.\smallskip

Hence, we see from example \ref{badcaseexample} that there is probably no versal space for non-trivial CR bundles. Unfortunately, we are unable to prove such a statement, because one cannot exclude that a construction completely different from ours yields a versal space. This is a crucial point that we hope to address in the future. If our intuition is correct, this would not only mean that Kodaira-Spencer and Kuranishi standard tools are not adapted to situations where there is no complete finite-dimensional space, but also that the concept of versality is not the good one. Linked to this is the question of the uniqueness of our Kuranishi type moduli space. In the classical setting of complex structures, versality implies uniqueness as a germ. However, for proper CR submersions over the circle, we are only able to prove a weak version of uniqueness in Corollary \ref{uniquenesscorollary}. And even this weak property becomes false for examples as example \ref{badcaseexample}.\smallskip

Rigidity is treated in Theorem \ref{rigiditytheorem} but this is indeed a very simple result. The only rigid case is the trivial example \ref{rigidexample} of a CR bundle with rigid fiber. In the classical setting of complex structures, rigidity results were the first important results to be proved. For that reason, it may seem natural to prove such results for each new deformation problem. However, when the situation is more flexible, as it is with CR structures, rigidity is not a good concept, because it so hardly happens, cf. the analogous case of \cite{Hamilton2}. Here, we derive it from connectedness results (Theorem \ref{connectedtheorem}, Corollaries \ref{concorI} and \ref{concorII}), which are by far more interesting applications of Theorem \ref{main}. \smallskip

The article ends with a section of examples.

\section{Kuranishi's Theorem as the existence of a local foliation.}
\label{Kuranishi}

In this Section, we state a version of Kuranishi's Theorem which is suited for our purposes. Although it is very close to the statements of \cite{Ku2}, this slight reformulation will be crucial in the proof of our main results. Indeed, we interpret Kuranishi's Theorem as describing a local foliated structure of finite codimension of the set of complex operators.
\vspace{5pt}\\
Now, let $H_{0}$ be the subspace of $\Sigma (TX^{diff})$ consisting of the real parts of the holomorphic vector fields of $X_{0}$. As $X_0$ is compact, its automorphism group is a finite-dimensional complex Lie group hence $H_0$ is finite-dimensional. 
Choose a decomposition
\begin{equation}
\label{splitting}
\Sigma (TX^{diff})= H_{0}\oplus L_{0}
\end{equation}
for some {\it closed} subspace $L_{0}$. Observe that such a closed complementary subspace always exists, since $H_0$ is finite-dimensional hence closed and since we are in a Hilbert space. Nevertheless, we want to emphasize that {\it we do not ask }\eqref{splitting} {\it to be orthogonal with respect to any product}. Indeed, in the sequel, we will use the fact that $L_0$ is not unique and that we can choose it.

\begin{remark}
It is important to keep in mind that the vector spaces of \eqref{splitting} are indeed $\mathbb C$-vector spaces through the identification between $\Sigma (TX^{diff})$ and $\Sigma (T)$; see also remark \ref{classicalremark}.
\end{remark}

Kuranishi's Theorem may be rephrased as:

\begin{theorem}
\label{KurTheorem}
There exists an open neighborhood $U_0$ of $J_0$ in $\mathcal E$, an open neighborhood $W_1$ of $0$ in $L_{0}$ and an analytic retraction
\begin{equation}
\label{Kurretraction}
\begin{CD}
U_0\cap\mathcal I @ > \Xi_0 > a.c. > K_0
\end{CD}
\end{equation}
such that
\begin{enumerate}
\item[(i)] The set $K_0$
is a (finite-dimensional) analytic set of (embedding) dimension at $0$ equal to 
$$
h^1(0):=\dim H^1(X_{0},\Theta)
$$
where $H^1(X_0,\Theta)$ denotes the first cohomology group with values in the sheaf $\Theta$ of germs of holomorphic vector fields.

\item[(ii)] The map
\begin{equation}
\label{Kuriso}
\begin{CD}
(J,\xi)\in K_0\times W_1 @ >\Phi_0>> e(\xi)\cdot J\in U_0\cap \mathcal I
\end{CD}
\end{equation}
is an a.c. isomorphism whose inverse has component in $K_0$ equal to $\Xi_0$.
\end{enumerate}
\end{theorem}

\begin{remark}
Isomorphism \eqref{Kuriso} must be thought of as a foliated chart. In other words, Kuranishi's Theorem states the existence of a finite codimension local foliated structure at each point of $\mathcal I$. Observe that it depends on the choice of $L$.
\end{remark}

\begin{remark}
\label{classicalremark}
In the classical presentation of \cite{Ku2}, a slightly different splitting  is used. From the decomposition \eqref{T}, Kuranishi defines $A^p$ as the space of $(0,p)$-forms of $X_{0}$ with values in $T$. Hence $A^0$ is the space of $(1,0)$-vectors and we may write
\begin{equation}
\label{classicalsplitting}
A^0=H^0\oplus ^{\perp} \kern-5pt H^0
\end{equation}
for $H^0$ the subspace of holomorphic vector fields and $^\perp H^0$ its orthogonal with respect to the $L^2$-norm for some fixed hermitian metric of $X_{0}$. And he encodes the small diffeomorphisms of $X^{diff}$ through the map
$$
\xi\in A^0\longmapsto e(\xi+\bar\xi).
$$
The version we present can easily be deduced from the classical one, the decomposition \eqref{splitting} playing the role of \eqref{classicalsplitting}. However, the crucial point is that \eqref{splitting} gives a splitting of the space of smooth vector fields, which is obviously independent of the complex structure $J_0$, whereas \eqref{classicalsplitting} gives a splitting of $A^0$, whose definition depends on $J_0$. So using \eqref{classicalsplitting} instead of \eqref{splitting} will allow us to compare different splittings based at different points. 
\end{remark} 

Let us give the core of the proof and explain why it is possible to replace the splitting \eqref{classicalsplitting} with the splitting \eqref{splitting}.
\vspace{5pt}\\
Let $\delta$ be the adjoint (from differential operator theory) with respect to a fixed hermitian metric on $X_0$ of the $\bar\partial$-operator extended to the forms with values in the holomorphic bundle $T$. This $\bar\partial$-operator acting on the spaces $A^p$ defines an elliptic complex, hence the associated Laplace-type operator $\Delta$ is elliptic. 
\vspace{5pt}\\
Also, by ellipticity, we have a direct sum decomposition
$$
A^1=\text{Im }\bar\partial\oplus\text{Ker }\delta.
$$
Using the splitting \eqref{classicalsplitting}, consider the smooth map
\begin{equation}
\label{mapKur}
(\xi,\omega)\in  ^{\perp} \kern-5ptH^0\times \text{Ker }\delta\longmapsto e(\xi+\bar\xi)\cdot \omega\in A^1.
\end{equation}
\begin{remark}
In order to make the previous decompositions precise, one should add that $A^1$ means the set of $1$-forms of class $W_2^r$, and $A^0$ means the $0$-forms of class $W_2^{r+1}$.
\end{remark}

A direct computation shows that its differential at $0$ is given by
\begin{equation}
\label{diffKur}
(\eta,\alpha)\in  ^{\perp} \kern-5ptH^0\times \text{Ker }\delta\longmapsto \alpha+\bar\partial\eta\in A^1
\end{equation}
and that \eqref{diffKur} is invertible with inverse given by
\begin{equation}
\label{Kurinverse}
\omega=\bar\partial \eta+\alpha\in\text{Im }\bar\partial\oplus\text{Ker }\delta\longmapsto
(G\delta\bar\partial\eta,\alpha)\in  ^{\perp} \kern-5ptH^0\times \text{Ker }\delta.
\end{equation}
Here $G$ denotes the Green operator associated to $\Delta$.
\vspace{5pt}\\ 
By application of the inverse function theorem on Hilbert spaces, the map \eqref{mapKur} is a local diffeomorphism. 
\vspace{5pt}\\
The Kuranishi space of $X_0$ is then defined as 
$$
K_0=\{\omega\in A^1\mid \bar\partial\omega-[\omega,\omega]=\delta\omega=0\}.
$$

The restriction of \eqref{mapKur} to $^{\perp} \kern-3pt H^0\times K_0$ gives a map similar to \eqref{Kuriso}, which appears in the classical statement of Kuranishi's Theorem. 
\vspace{5pt}\\
If we use now the splitting \eqref{splitting}, we just have to modify the previous formulas as follows.
\vspace{5pt}\\
We consider the map
\begin{equation}
\label{newKurmap}
(\xi,\omega)\in  L_0\times \text{Ker }\delta\longmapsto e(\xi)\cdot \omega\in A^1
\end{equation}
instead of \eqref{mapKur}, whose differential at $0$ is
$$
(\eta,\alpha)\in  L_0\times \text{Ker }\delta\longmapsto \alpha+\bar\partial \tau\eta\in A^1.
$$
This is completely analogous to \eqref{diffKur}, the only difference being the use of the identification
$$
\tau : \xi\in \Sigma (TX^{diff})\longmapsto \xi-iJ_0\xi\in A^0.
$$
This identification maps $\Sigma (TX^{diff})$ onto $A^0$, and $H_0$ onto $H^0$. But it does not map $L_0$ onto $^{\perp} \kern-3pt H^0$. Nevertheless, since $\tau L_0$ and  $^{\perp} \kern-3pt H^0$ are complementary to the same finite-dimensional space, they are isomorphic, so that we can twist it into an identification
\begin{equation}
\label{tildetau}
\tilde\tau : H_0\oplus L_0\longrightarrow H^0\oplus ^{\perp} \kern-5ptH^0
\end{equation}
which preserves the direct sum decompositions. More precisely, we define $\tilde\tau$ as the map
\begin{equation}
\label{tildetau2}
\tilde\tau (\xi_0\oplus\xi_{L_0})=\tau\xi_0\oplus(\tau\xi)^\perp
\end{equation}
where
$$
\xi=\xi_0\oplus\xi_{L_0}\in H_0\oplus L_0\quad\text{ and }\quad
\tau\xi=(\tau\xi)^0\oplus (\tau\xi)^\perp\in H^0\oplus ^{\perp} \kern-5ptH^0.
$$
From \eqref{tildetau} and \eqref{tildetau2}, we infer that the formula for the inverse of \eqref{newKurmap}, analogous to \eqref{Kurinverse}, is
$$
\omega=\bar\partial \tau\eta+\alpha\in\text{Im }\bar\partial\oplus\text{Ker }\delta\longmapsto
(\tilde\tau^{-1}G\delta\bar\partial\tau\eta,\alpha)\in  L_0\times \text{Ker }\delta.
$$

This shows that the map \eqref{newKurmap} is a local diffeomorphism. Its restriction to $L_0\times K_0$ gives the map \eqref{Kuriso}. The first component of its inverse gives the map \eqref{Kurretraction}.

\begin{definition}
\label{Kurdomain}
We call a map $\Xi_0$ as defined in \eqref{Kurretraction} a {\it Kuranishi map} based at $J_0$. The open set $U_0$ is called a {\it Kuranishi domain}. The analytic space $K_0$ is called the {\it Kuranishi space} of $X_{0}$. 
\end{definition}

We note the following uniqueness property.

\begin{corollary}
The Kuranishi space $K_0$ of $X_{0}$ is unique in the following sense.
\begin{enumerate}
\item[(i)] If $L'$ is another closed complementary subspace to $H_{0}$, then the corresponding space $K'_0$ is a.c. isomorphic to $K_0$.

\item[(ii)] If $U'$ is another neighborhood of $J_0$ in $\mathcal E$, then the restrictions of $\Xi_0$ and of the corresponding map $\Xi'_0$ to $U\cap U'\cap\mathcal I$ have a.c. isomorphic images.
\end{enumerate}
\end{corollary}
In particular, it is unique as a germ of analytic space at $0$. In this paper, thinking of this corollary, we say that $K_0$ is {\it the} Kuranishi space of $X_0$, even if, strictly speaking, it depends on $U_0$ and on $L_0$.
\vspace{5pt}\\
In the sequel, when we will refer to Kuranishi's Theorem, we will always refer to the version given in Theorem \ref{KurTheorem}.

\section{Foliation of a neighborhood of a compact set of $\mathcal I$.}
\label{foliationcompact}

Recall that $\text{Diff}^0 (X^{diff})$ acts on $\mathcal E$. Given a closed vector subspace $L$ of the space $\Sigma (TX^{diff})$, 
we will see that, under some conditions, it induces a foliation of certain open sets of $\mathcal I$ whose leaves are included in the orbits of the $\text{Diff}^0(X^{diff})$ action. We first need the following technical result.

\begin{lemma}
\label{homotopiclemma}
Let $J$ and $J'$ be two homotopic complex structures. Then $\Sigma (T_J)$ and $\Sigma (T_{J'})$ are isomorphic $\mathbb C$-vector spaces.
\end{lemma}

\begin{proof}
We may assume that $T_{J'}$ is represented by the graph of a $1$-form from $T_J$ to $\overline T_J$; otherwise just choose a finite number of points  $J_i$ on the path from $J$ to $J'$ and repeat the argument inductively. The projection on $T_J$ parallel to $\overline T_J$ identifies $T_{J'}$ and $T_J$ as $\mathbb C$-vector bundles over $X^{diff}$. Hence $\Sigma (T_J)$ and $\Sigma (T_{J'})$ are isomorphic $\mathbb C$-vector spaces.
\end{proof}

Thanks to this Lemma, the structure of complex vector space induced on $L$ through \eqref{splitting} and the identification between $TX^{diff}$ and $\overline T_J$ does not depend on $J$, provided we stay in a single connected component of complex operators.

\begin{theorem}
\label{foliation}
Let $C$ be a compact set of $\mathcal E$. 
Then, for any sufficiently small neighborhood $U$ of $C$ in $\mathcal E$, there exists a closed subspace $L$ of $\Sigma (TX^{diff})$ of finite codimension such 
that $L$ induces a foliation of $U\cap\mathcal I$.
More precisely, there exists a finite open covering
\begin{equation}
\label{Ucovering}
U=U_1\cup \hdots \cup U_k
\end{equation}
such that
\begin{enumerate}
\item[(i)] For every $i$, there exists an analytic space $K_i$ such that $U_i\cap\mathcal I$ is a.c. isomorphic to the product of some fixed neighborhood $W$ of $0$ in $L$ with $K_i$. 

\item[(ii)] When defined, the composition of two such isomorphisms preserves the plaques $W_i\times\{Cst.\}$.
\end{enumerate}
\end{theorem}

\begin{remark}
Note that the fact that the leaves are included in the orbits of the $\text{Diff}^0(X^{diff})$ action comes from the fact that the above isomorphisms are a.c. Also note that the leaves are modeled onto $L$, which is a well-defined $\mathbb C$ vector space by Lemma \ref{homotopiclemma}.
\end{remark}

\begin{definition}
\label{Lfoliationdef}
In the setting of Theorem {\rm \ref{foliation}}, we say that $L$ {\it foliates} $U\cap\mathcal I$ and we call the associated foliation the $L$-{\it foliation} of $U\cap\mathcal I$.
\end{definition}

\begin{remark}
The $L$-foliation is transversely modeled on $K_i$, an analytic space which depends on the chart $U_i$. So strictly speaking, we should talk of a lamination rather than a foliation.
\end{remark}

\begin{proof}
Begin with choosing an open neighborhood $U$ of $C$ in $\mathcal E$ and a finite open covering of $U$ by Kuranishi domains $U_1,\hdots , U_k$ based at $J_1,\hdots ,J_k$, points of $C$. Set
$$
H_i:=H_{J_i}
$$
and denote by $\text{Kur}_i$ the Kuranishi space of $J_i$. Assume that we used the classical orthogonal splitting, that is, in the splitting \eqref{classicalsplitting}, we took $L_i$ as the orthogonal complement to $H_i$. It follows from Kuranishi's Theorem that for every $i$ and every $J\in U_i\cap\mathcal I$, we have
\begin{equation}
\label{LiHj}
L_i\cap H_J=\{0\}.
\end{equation}
Let $\mathcal S$ be the set of closed vector subspaces $S$ of $\Sigma (TX^{diff})$ having finite codimension and satisfying 
\begin{equation}
\label{S}
S\cap H_J=\{0\}\qquad\text{for all }J\in U.
\end{equation}
We claim that, shrinking $U$ if necessary, we may assume that $\mathcal S$ is not empty. Indeed, let 
$$
S:=(H_1+\hdots+ H_k)^\perp .
$$
Then $S$ is orthogonal to each $H_i$, hence it is included in each $L_i$. Property \eqref{LiHj} implies that $S$ belongs to $\mathcal S$. As a consequence, we may 
choose $L\in \mathcal S$ having minimal codimension. We claim that this $L$ satisfies the requirements of Theorem \ref{foliation}. 
\vspace{5pt}\\
Since $L$ has finite codimension in $\Sigma (TX^{diff})$, we may choose some finite-dimen\-sional vector subspaces $\tilde H_i$ such that, for all $i$, we have
$$
\Sigma (TX^{diff}) =L\oplus \tilde H_i\oplus H_i .
$$
Set
\begin{equation}
\label{tildeLi}
\tilde L_i:=L\oplus \tilde H_i.
\end{equation}
We may then replace $L_i$ with $\tilde L_i$ and obtain new Kuranishi maps and domains based at $J_i$. If they do not cover $C$, we just add a finite number of 
extras Kuranishi maps. This is possible thanks to \eqref{S}. To simplify notations, we still denote by the same symbols $U_i$ this refined covering. 
Remembering \eqref{tildeLi}, Kuranishi's Theorem implies that these Kuranishi maps induce a.c. isomorphisms between $U_i\cap \mathcal I$ and the product of the 
finite-dimensional analytic space 
$$
K_i:=\text{Kur}_i\times B_i
$$ 
(for $B_i$ an open neighborhood of $0$ in $\tilde H_i$) with $W$, an open neighborhood of $0$ in $L$ that we may assume to be the same for all $i$. Hence they give foliated charts as wanted in (i). 
\vspace{5pt}\\ 
Observe that the case where at least one of the Kur$_i$ is not reduced is treated exactly in the same way.
\vspace{5pt}\\
Moreover, still by Kuranishi's Theorem, the plaques $W\times\{Cst\}$ correspond exactly to the local equivalence classes of the relation
$$
J\sim J'\iff J'=e(\xi)\cdot J\quad\text{ for some }\xi\in W
$$ 
hence have an intrinsic geometric meaning and are preserved by the changes of charts. 
This proves (ii).
\end{proof}

\begin{definition}
Let $U\cap \mathcal I$ be $L$-foliated. Then {\it a foliated atlas} of $U\cap\mathcal I$ is a finite covering \eqref{Ucovering} satisfying the conclusions of Theorem \ref{foliation}. It is an {\it adapted foliated atlas} if it is a foliated atlas and an adapted covering in the sense of definition \ref{adapteddef}.
\end{definition}
 
Observe that $L$ depends on $U$. Also, observe that, given such an $U$, there are uncountably many isomorphic $L$ such that $L$ foliates $U\cap\mathcal I$ (just take one and perturb it slightly). Unfortunately, we do not know if two such foliations are equivalent. But we have the following uniqueness property on $L$.

\begin{proposition}
\label{Luniqueness}
The subspace $L$ is unique in the following sense. If $L$ and $L'$ are two vector subspaces of $\Sigma (TX^{diff})$ such that
\begin{enumerate}
\item[(i)] $L$ and $L'$ foliate $U\cap\mathcal I$.

\item[(ii)] Both $L$ and $L'$ are elements of $\mathcal S$ of minimal codimension.
\end{enumerate}
Then $L$ and $L'$ are isomorphic.
\end{proposition}

Observe that, if $L$ foliates $U\cap\mathcal I$, then $L$ must be an element of $\mathcal S$. And of course, if they do not have the same codimension, we cannot expect the foliations to be isomorphic. So condition (ii) in the statement of Proposition \ref{Luniqueness} is not a restriction.

\begin{proof}
This is standard linear algebra. Set
$$
I=L\cap H \qquad J=L'\cap H
$$
where 
$H:=L^\perp + (L')^\perp$. Observe that both $L+H$ and $L'+H$ generates $\Sigma(TX^{diff})$. 
Then choose $L_0$ and $L'_0$ closed such that
$$
\begin{aligned}
\Sigma (TX^{diff})=&H\oplus L_0\qquad\qquad\text{ with } L=L_0\oplus I\cr
=&H\oplus L'_0\qquad\qquad\text{ with } L'=L'_0\oplus J.
\end{aligned}
$$ 
So $L_0$ and $L'_0$ are isomorphic as complementary subspaces of the same subspace $H$, and have thus the same codimension. So $I$ and $J$ are isomorphic and we may extend the isomorphism between $L_0$ and $L'_0$ to an isomorphism between $L$ and $L'$.
\end{proof}

\section{Foliation of a neighborhood of a path.}
\label{foliationpath}

Let $c$ be a continuous path into $\mathcal I$ satisfying Hypothesis \ref{noselfint} and let $U$ be a connected neighborhood of $c$ in $\mathcal E$ which is $L$-foliated for some $L$. In this particular case, we can give a much more precise description of the $L$-foliation.

\begin{theorem}
\label{foliationtrivial}
Let $L$ foliating $U\cap \mathcal I$. Then, if $U$ is small enough, the foliation is given by the level sets of an analytic morphism. To be more precise, there exists an analytic space $K_U$ and an a.c. morphism
$$
\begin{CD}
U\cap\mathcal I@ >\Xi_U >a.c. >K_U
\end{CD}
$$
such that 
\begin{enumerate}
\item[(i)] $K_U$ is the leaf space of the $L$-foliation and the leaves are given by the level sets of $\Xi_U$.

\item[(ii)] The map $\Xi_U$ is locally a projection: in the neighborhood $V$ of any point $x\in U\cap\mathcal I$, we have a commutative diagram
\begin{equation}
\label{triviality}
\begin{CD} 
x\in V@ >\text{a.c. isomorphism}>> \Xi_U(V)\times W\cr
@ V\Xi_U VV @ VV \text{1st projection}V\cr
\Xi_U(V)\subset K_U@ >>\text{Identity}> \Xi_U(V)\subset K_U
\end{CD}
\end{equation}
for some open neighborhood $W$ of $0$ in $L$.
\vspace{5pt}\\
Moreover, the $L$-foliation is smoothly trivial, that is there exists a smooth injection
\begin{equation}
\label{trivialinjection}
\begin{CD}
K_U@ >i_U>a.c. > U\cap\mathcal I
\end{CD}
\end{equation}
and, up to shrinking $U$ if necessary, a diffeomorphism
\begin{equation}
\label{trivialdiff}
\begin{CD}
(J,f)\in K_U\times W_U@ > \Phi_U >a.c. > f\cdot i_U(J)\in U\cap\mathcal I
\end{CD}
\end{equation}
where $W_U$ is a subset of $\text{\rm Diff}^0(X^{diff})$ containing the identity.
\end{enumerate}
\end{theorem}

In other words, the $L$-foliation is smoothly trivial, hence its leaf space is an analytic space, cf. section \ref{outline}.
\begin{remark}
The open set $U$ may contain "large" open sets of leaves, so that we cannot ensure that they can be completely encoded via the map $e$. This is typically the case if the path $c$ is included in a single leaf. This explains why the set $W_U$ of \eqref{trivialdiff} is included in $\text{Diff}^0(X^{diff})$ and not in $L$. 
\end{remark}

\section{Proof of Theorem \ref{foliationtrivial}.}
\label{prooftrivial}

Recall the line of arguments presented in section \ref{outline}.\vspace{5pt}\\
If $c$ is constant or is contained in a single Kuranishi domain, then there is nothing to do: we just take $U$ as a Kuranishi domain based at $c(0)$ and the map \eqref{Kuriso} given by Kuranishi's Theorem as trivialization chart  \eqref{triviality}. This map is indeed global, so that \eqref{trivialdiff} is satisfied with an analytic isomorphism.\vspace{5pt}\\
So assume that $c$ is not contained in a single Kuranishi domain.\vspace{5pt}\\
We just make use of Theorem \ref{foliation}. We first assume that $c$ is not a loop. We cover it by an adapted foliated covering. 
Furthermore, we assume that the $U_i$ are sufficiently small so that, for any pair $(x,y)$ of points of $U_i\cup U_{j}$ with $U_i\cap U_j\not =\emptyset$, if $x$ and $y$ belong to the same orbit in $U_i\cup U_{j}$, we have
\begin{equation}
\label{furtherass}
y=e(\xi)\cdot x\qquad\text{for some }\xi\in W_i\cap W_{j}.
\end{equation}
We want to construct the leaf space $K_U$ by gluing each $K_i$ to the next one using the changes of foliated charts. The difficult point is to ensure it is Hausdorff. To do that, we proceed by induction.\vspace{5pt}\\
Firstly, consider, for $i=1,2$, the two local charts 
\begin{equation}
\label{localcharts}
 \begin{CD}
 \Phi_i\ :\ K_i\times W_i@ >a.c.>> U_i\cap\mathcal I
\end{CD}
\end{equation}
and the associated change of chart
$$
\begin{CD}
(z,v)@ > \Phi_{12}:=\Phi_2^{-1}\circ\Phi_1> a.c. > (\Psi (z),\chi (z,v)).
\end{CD}
$$
By abuse of notation, we will denote by the same symbol $K_i$ and the image $\Phi_i(K_i\times\{0\})$. We assume, composing $\Phi_1$ and $\Phi_2$ with translations on the $W$-factor if necessary, that $K_1$ and $K_2$ are disjoint in $\mathcal I$.\vspace{5pt}\\
Consider the space obtained by gluing $K_1$ to $K_2$ through the map $\Psi$. We claim that, up to shrinking $U_1$ and $U_2$, it is an analytic space. Hence the foliation restricted to $U_1\cup U_2$ satisfies the properties of Theorem \ref{foliationtrivial} and the first step of the induction is done. This can be proved as follows.\vspace{5pt}\\
We define the open set $V_1\subset K_1$ by the following relation
\begin{equation}
\label{relation}
x\in V_1\iff \exists v\in W_1\quad\text{ such that } \Phi_1(x,v)\in U_2\cap\mathcal I.
\end{equation}
By \eqref{furtherass}, observe that $V_1$ contains all the points whose leaf in $U_1$ meets $U_{2}$. 
\begin{remark}
Recall that the leaf in $U_1$ of a point $J$ is the {\it connected component} of $J$ in the intersection of $U_1$ with the leaf passing through $J$.
\end{remark}
Denoting as usual by $\Xi_1$ the map contracting $U_1\cap\mathcal I$ onto $K_1$, observe that $\Xi_1(c\cap U_1)\cap V_1$ is non-empty and connected since it is equal to $\Xi_1(c\cap U_1\cap U_2)$ (recall Hypothesis \ref{noselfint}). Then, define
$$
\begin{CD}
x\in V_1@ >g>a.c.> \Xi_2\circ \Phi_1(x,v)\in K_2
\end{CD}
$$
where $v$ is any vector of $W_1$ such that \eqref{relation} holds. The previous map does not depend on the particular choice of $v$ because of \eqref{furtherass}.\vspace{5pt}\\
This is the gluing map we want to use. Set $V_2=g(V_1)\subset K_2$.
It follows from the proof of Theorem \ref{foliation} that 
there exists an analytic map
$$
\eta \ :\ V_1\subset K_1\longrightarrow L
$$
such that 
\begin{equation}
\label{geta}
g(J)= e(\eta(J))\cdot J.
\end{equation}
Now, choose some open set $V'_1$  included in $V_1$ such that  $\Xi_1(c\cap U_1)\cap V'_1$ is still non-empty and connected. Let $\chi$ be a bump function defined on $U_1$, equal to $1$ on $V'_1$ and to $0$ on $K_1\setminus V_1$. Let $V'_2$ be the image $g(V'_1)$. Define
$$
\eta_1:=\chi\cdot \eta\ :\ U_1\longrightarrow L.
$$
Then the map $g_1:=e(\eta_1)$ is an a.c. diffeomorphism from $K_1$ to its image. Moreover, we may assume that $g_1(K_1)$ intersects $K_2$ along
$$
\tilde V'_2=\{y\in \overline{V'_2}\mid L^{U_1\cup U_2}_y\cap \overline{V'_1}\not =\emptyset\}.
$$
where $L^{U_1\cup U_2}_y$ denotes the leaf of $y$ in $U_1\cup U_2$, and where the closure is taken in $K_2$ (respectively $K_1$).\vspace{5pt}\\
We would like that the set $\tilde V'_2$ is open and that the union $g_1(K_1)\cup K_2$ cuts every leaf of $U_1\cup U_2$ into a single point. Taking into account the definition of $\tilde V'_2$, this would imply that the image through $g_1$ of $\overline{V_1}'\setminus V'_1$ does not belong to $K_2$. From this, it is easy to check that $K_1\cup_{g_1}K_2$ would be Hausdorff, hence an analytic space since it is obtained by gluing two analytic spaces along an open set; also this would imply that it is the leaf space of $U_1\cup U_2$.
\begin{remark}
\label{clarifyremark}
To clarify the proof, let us emphasize the following obvious but crucial point: $K_1\cup_{g_1}K_2$ will be constructed as an {\it abstract} analytic space {\it homeomorphic} to $g_1(K_1)\cup K_2$. But obviously, this last set is not an analytic subspace of $\mathcal I$, cf. section \ref{outline}.
\end{remark}

\begin{figure}[H]
    \centering
    \includegraphics[width=.6\textwidth]{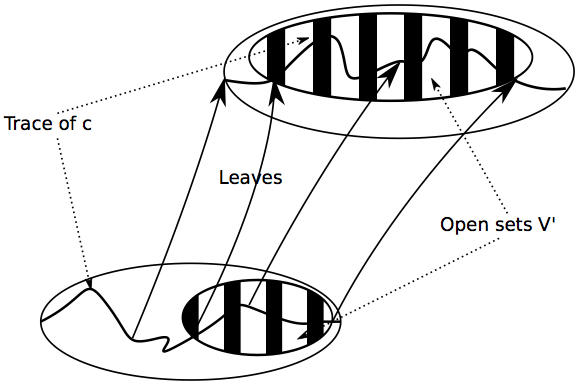}
    \caption{Non-Hausdorff gluing.}
    \label{NHaus}
\end{figure}

\begin{figure}[H]
    \centering
    \includegraphics[width=.6\textwidth]{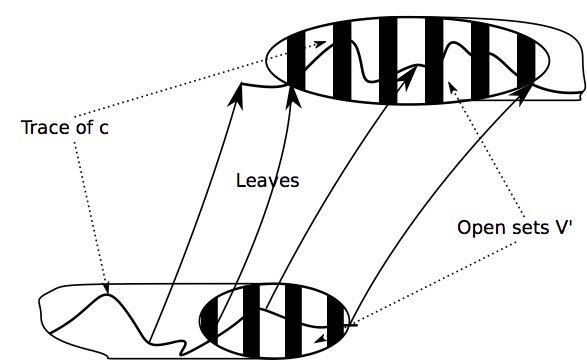}
    \caption{Hausdorff gluing.}
    \label{Haus}
\end{figure}

Nevertheless, it is not true in general. Indeed, it fails every time that $\tilde V'_2$ contains a point of $\overline{V'_2}\setminus V'_2$.\vspace{5pt}\\
We claim that it is enough to shrink $U_1$ and $U_2$ to ensure
\begin{equation}
\label{tilde2}
\tilde V'_2=V'_2
\end{equation}
and thus to solve our problem.\vspace{5pt}\\
In Figure \ref{NHaus}, the big ellipses represent $K_1$ and $K_2$ and the small ones represent $V'_1$ and $V'_2$. The arrows on the leaves just suggest the identification. The glued space is obtained by gluing along the {\it open} shaded parts. Clearly, since $V'_1$ has boundary points in $K_1$ which correspond through the identification to boundary points of $V'_2$ in $K_2$, the resulting space is not Hausdorff.\vspace{5pt}\\
The claim is that it is possible to shrink $U_1$ and $U_2$ so that the picture is like Figure \ref{Haus}. The big ellipses represent the domains of $K_1$ and $K_2$ onto which the new sets $U_1$ and $U_2$ retract. In this new case, the boundary points of $V'_1$ are not in correspondance with the boundary points of $V'_2$, hence the gluing is Hausdorff. Notice that, in the case of $U_2$ (the upper ellipse), a subset of the former trace of $c$ in $K_2$ is now out of the domain. This has no consequence since this part of $c$ has still a trace in $K_1$, but it helps reducing the boundary of $V'_2$.
First, observe that shrinking $U_i$ to, say, $U'_i$, we may assume that it transforms \eqref{localcharts} into
\begin{equation}
\label{newlocalcharts}
\begin{CD}
\Phi_i \ :\ K'_i\times W'_i\subset K_i\times W_i@ > a.c. >> U'_i\cap\mathcal I
\end{CD}
\end{equation}
However, $W'_i\subset W_i$ {\it may not contain} $0$. Then $g_1(K'_1)$ intersects $K'_2$ along $\tilde V'_2\cap K'_2$. The claim is that, after shrinking and after taking the intersection with $K'_2$, we may assume that \eqref{tilde2} holds.
To go further, we distinguish cases.\vspace{5pt}\\
$\underline{\text{1st case}}$.
We assume that $\Xi_i (c\cap U_i)$ does not intersect the boundary $\overline{V'_i}\setminus V'_i$. This is typically the case that $c$ is included in a single leaf of $L$, hence its trace in each $K_i$ is a single point belonging to $V'_i$. Then there is no problem. We set
$$
U_i=\Xi_i^{-1} (V'_i)\quad\text{ for }\quad i=1,2.
$$
Geometrically, $K_1$ and $K_2$ are simply identified one to the other.\vspace{5pt}\\
$\underline{\text{2nd case}}$.
We assume that the intersection of $\Xi_i (c\cap U_i)$ with the boundary $\overline{V'_i}\setminus V'_i$ is connected. Then it is enough to shrink $U_2$ setting
$$
U_2=\Xi_2^{-1} (V'_2).
$$
Geometrically, $K_2$ is identified to an open subset of $K_1$.\vspace{5pt}\\
$\underline{\text{3rd case}}$.
We assume that the intersection of $\Xi_i (c\cap U_i)$ with the boundary $\overline{V'_i}\setminus V'_i$ is disconnected. This is the case treated in Figures \ref{NHaus} and \ref{Haus}.\vspace{5pt}\\
Set
$$
I=\{t\in [0,1]\mid \Xi_2(c(t))\cap \overline{V'_2}\setminus V'_2\not =\emptyset\}.
$$
It has at least two connected components and we may find two disjoint open intervals $I_1$ and $I_2$ such that
\begin{equation}
\label{I12}
I\subset I_1\sqcup I_2.
\end{equation}
Then, shrink $U_1$ and $U_2$ such that Hypothesis \ref{noselfint} and \eqref{furtherass} are still satisfied, as well as the additional hypothesis
\begin{equation}
\label{U'}
U'_i\cap c(I)\subset c(I_i).
\end{equation}
Properties \eqref{I12} and \eqref{U'} imply that \eqref{tilde2} holds
after taking the intersection with $K'_2$. \vspace{5pt}\\
In other words, set
$$
K_{12}:=K'_1\cup_{g_1} K'_2\quad\text{ with }g_1\ :\ K'_1\cap V'_1\to K'_2\cap V'_2.
$$
Then $K_{12}$ is an analytic space. It is by construction the leaf space of the $L$-foliation restricted to $U'_1\cup U'_2$.

\begin{remark}
\label{notreducedremark}
If $K_1$ (or equivalently $K_2$)  is not reduced, then we perform the previous construction with their reduction and put on each component of the resulting space the common multiplicity of the corresponding components of $K_1$ and $K_2$.
\end{remark}

Besides the two maps
$$
\begin{CD}
U'_1\cap\mathcal I@ >\Phi_1^{-1} > a.c. > K'_1\times W'_1@ >\text{1st projection} >> K'_1
\end{CD}
$$
and
$$
\begin{CD}
U'_2\cap\mathcal I@ >\Phi_2^{-1} > a.c. > K'_2\times W'_2@ >\text{1st projection} >> K'_2
\end{CD}
$$
glue into a single map
\begin{equation}
\label{singlechart}
\begin{CD}
\Xi :\ (U'_1\cup U'_2)\cap\mathcal I@ > a.c. >> K_{12}
\end{CD}
\end{equation}
with charts \eqref{trivialinjection} by construction.\vspace{5pt}\\
It follows now from \eqref{singlechart} that the inclusions
$$
i_1\ :\ z\in K'_1 \longmapsto\Phi_1(z,\eta_1(z))\in U'_1\cap\mathcal I
$$
and 
$$
i_2\ :\ z\in K'_2 \longmapsto\Phi_2(z,0)\in U'_2\cap\mathcal I.
$$
glue naturally into a smooth inclusion
$$
\begin{CD}
i_{12}\ :\ K_{12}@ > a.c. >> (U'_1\cup U'_2)\cap\mathcal I
\end{CD}
$$
yielding \eqref{trivialinjection} and a smooth trivialization \eqref{trivialdiff} (shrinking $U$ if necessary).\vspace{5pt}\\
Repeating the process, we construct the map $\Xi_U$ as desired as well as the smooth trivialization. So we are done.\vspace{5pt}\\
Assume now that $c$ is a loop. We use once again a foliated adapted covering and assume \eqref{furtherass}.
We proceed as before, but we have now to perform an ultimate gluing between $U_k$ and $U_1$ to obtain $K_U$, still using the changes of charts of Theorem \ref{foliation}. Let $K$ denote the analytic space obtained by making all the gluings except for the last one. We also have a smooth  map 
$$
\begin{CD}
i\ :\ K@ > a.c. >> U\cap \mathcal I
\end{CD}
$$
analogous to \eqref{trivialinjection}. The last gluing to perform is defined through an analytic map
$$
\begin{CD}
V_k\subset K_k@ >g>a.c.>K_1
\end{CD}
$$
which is equal as before to $e(\eta)$ for some analytic map $\eta$ from $K_k\cap U_1$ into $L$ because of \eqref{furtherass}. So we may proceed as before and extend smoothly $\eta$ to
$$
\begin{CD}
\eta_1\ :\ K@ >>> L
\end{CD}
$$
equal to $\eta$ on $V'_k\subset K_k$ for some open set $V'_k$ included in $V_k$. We assume it meets $\Xi_k(c\cap U_k)$.  Now, up to shrinking the $U_i$'s, we have that
$e(\eta)\cdot K$ is a smooth global transverse section to the $L$-foliation on $U\cap\mathcal I$. This proves at the same time that the space $K_U$ obtained from $K$ after performing the last gluing is homeomorphic to $e(\eta)\cdot K$, hence Hausdorff and by construction an analytic space; and that the foliation is smoothly trivial.\vspace{5pt}\\
The injection \eqref{trivialinjection} and the trivialization \eqref{trivialdiff} are then obtained as before.\vspace{5pt}\\
We note the following uniqueness property
\begin{corollary}
\label{KUuniqueness}
The analytic space $K_U$ is unique up to a.c. isomorphism, that is does not depend on the choice of the adapted covering and of the Kuranishi maps.
\end{corollary}

\begin{proof}
This is a direct consequence of the fact that $K_U$ is the leaf space of the $L$-foliation  restricted to $U$. Hence it is unique.
\end{proof}

For the same reason, we also have

\begin{corollary}
\label{KUuniqueness2}
Let $U$ and $U'$ be two connected neighborhoods of $c$ for which a smooth trivialization {\rm \eqref{trivialdiff}} exists. Then the restrictions of $K_U$ and $K_U'$ to $U\cap U'\cap \mathcal I$ (via the trivializations {\rm \eqref{trivialdiff}}) are a.c. isomorphic.
\end{corollary}

However, it is worth to emphasize that $K_U$ {\it depends on the choice of} $L$.

\section{The Kuranishi type moduli space of a proper CR submersion.}
\label{Kuranishitype}

We are now in position to prove the main result of this paper: a statement analogous to Kuranishi's Theorem for $\mathcal I(\mathcal X^{diff})$.\vspace{5pt}\\
Let $\mathcal X_0$ be a CR submersion compatible with $\pi^{diff}$, represented by an element $c_0$ of $\mathcal I(\mathcal X^{diff})$. Assume Hypothesis \ref{noselfint}. Identify a neighborhood of $c_0$ in $\mathcal E(\mathcal X^{diff})$ with $C^\infty (\mathbb S^1, U_\phi)$ as explained in subsection \ref{CRsubss}. Choose a closed vector subspace $L$ of $\Sigma (TX^{diff})$ satisfying \eqref{S} for all $J$ in the image of $c_0$ and having minimal codimension for this property. 
We have (compare with Theorem \ref{KurTheorem})

\begin{theorem}
\label{main}
Shrinking $U_\phi$ if necessary, we can find a finite-dimensional analytic space $K_{c_0}$ and an analytic map
\begin{equation}
\label{pathretraction}
\begin{CD}
U_\phi\cap\mathcal I @ > \Xi_\phi > a.c. > K_{c_0}
\end{CD}
\end{equation}
such that
\begin{enumerate}
\item[(i)] The (embedding) dimension at $c_0(t)$ of the space
$K_{c_0}$
is equal to 
$$
\begin{aligned}
&h_1(t)+\text{\rm codim }L-h_0(t)+1\quad\text{ if }\mathcal X_0 \text{ is a non-trivial CR bundle}\cr
&h_1(t)+\text{\rm codim }L-h_0(t)\qquad\text{ otherwise} 
\end{aligned}
$$
where $\text{\rm codim }L$ is 
the codimension of $L$ in $\Sigma (TX^{diff})$.

\item[(ii)] If we are not in the special case of (i), the analytic set $K_{c_0}$ is the leaf space of the $L$-foliation of $U_\phi$.  Otherwise, there exists a closed subspace $L'$ of $\Sigma (TX^{diff})$ contained in $L$ as a codimension-one subspace and such that $K_{c_0}$ is the leaf space of the $L'$-foliation of $U_\phi$.
\end{enumerate}
\end{theorem}

As a consequence, the (infinite-dimensional) analytic space $C^\infty (\mathbb S^1, K_{c_0})$ plays the role of the Kuranishi space $K_0$ in the classical case. Hence we define

\begin{definition}
The loop space $C^\infty (\mathbb S^1, K_{c_0})$ is called a {\it Kuranishi type moduli space} of $\mathcal X_0$. We denote it by $K^g$.
\end{definition}

\begin{proof}
We keep the same notations as in the previous sections and recall that $U_\phi$ is defined as the quotient of some open neighborhood $U$ of $c_0$ by $\phi$. By Theorem \ref{foliationtrivial}, shrinking $U$ and thus $U_\phi$ if necessary, attached to $U$ is an analytic space $K_U$ together with an analytic a.c. map of $U$ onto $K_U$ such that $K_U$ is the leaf space of the $L$-foliation of $U$. We want now to define an analytic space $K_{c_0}$ attached to $U_\phi$. 
\vspace{5pt}\\
If $\phi$ is the identity, there is nothing to do. We have $U_\phi=U$ and we take $K_{c_0}=K_U$, that is we take exactly the analytic space given by Theorem \ref{foliationtrivial}.\vspace{5pt}\\
To do the general case, it would be natural to define $K_{c_0}$ as the quotient of $K_U$ by the action of $\phi$. However, the resulting quotient space is not always an analytic space and we have to consider two different cases.
\vspace{5pt}\\
Indeed, in Hypothesis \ref{noselfint}, we imposed the condition of $\phi$ not being a biholomorphism. This forces $c_0$ to have distinct endpoints. It follows that, when gluing  $U_k$ to $U_1$ through $\phi$, the resulting quotient space is Hausdorff (at least after shrinking). The fact that $U_k$ and $U_1$ may be supposed to be disjoint is fundamental in this process. In the same way, when performing the same gluing onto $K_U$, we must ensure that the glued pieces corresponding to $K_1$ and $K_k$ in $K_U$ are disjoint to obtain an analytic space. This is possible (shrinking $U$ if necessary) if and only if the image $c_U$ of $c_0$ in $K_U$ is not a loop.
\vspace{5pt}\\
Indeed, we may assume that
$$
K_k=\phi\cdot K_1
$$
and, after identification between open sets of $K_k$ and $K_1$ and open sets of $K_U$, this induces a well-defined analytic a.c. isomorphism between two open sets of $K_U$. As we just told, if we may assume that, after shrinking, these two pieces are disjoint, then we may proceed as in the proof of Lemma \ref{loopstructureIIlemma}  (construction of $U_\phi$) and ensure that the gluing occurs exclusively on these open sets and that the resulting analytic space is the desired leaf space. Hence we are done.
\vspace{5pt}\\
Assume that $\mathcal X_0$ is not a CR bundle. Then, by Fischer-Grauert Theorem, we can find $t\not = t'$ such that $X_{c_0(t)}$ is not biholomorphic to $X_{c_0(t')}$. This implies that $c_U$ is not a constant path. It may of course be a loop, which is exactly the situation we would like to avoid, but since it is not a constant loop, we claim that, shrinking $U$, it becomes a path. Indeed, assume that $c_U$ is a loop. Then this means not only that $c_0(0)$ and $c_0(1)$ are in the same leaf of $\mathcal I$, but also that they are in {\it the same leaf of }$U$. Either $c_0(t)$ or $c_0(t')$ must be in a different $L$-orbit than that of the two endpoints. Say it is $c_0(t)$. By shrinking $U$, we may assume that $c_0(t)$ belongs to some domain $U_p$ with $p$ different from $1$ and from $k$, and that this $U_p$ does not intersect the leaf of $c_0(0)$. Hence the intersection of this orbit with $U$ is disconnected and $c_0(0)$ from the one hand, and $c_0(1)$ from the other hand, belong to two different connected components. In other words, the common leaf of $c_0(0)$ and $c_0(1)$ in $\mathcal I$ disconnects into (at least) two leaves in $U$, one passing through $c_0(0)$, and the other passing through $c_0(1)$. Because of the trivialization \eqref{localcharts}, this prevents their images $c_U(0)$ and $c_U(1)$ to be the same point of $K_U$. So in this case, we may define $K_{c_0}$ as the quotient of $K_U$ by the action of $\phi$.
\vspace{5pt}\\
Assume now that $\mathcal X_0$ is a CR bundle. Then, we cannot exclude that $c_U$ is the constant loop even after shrinking $U$ (cf. Example \ref{badcaseexample}). The quotient of $K_U$ by the action of $\phi$ occurs in the neighborhood of the {\it point} $c_U$, which is fixed by $\phi$. As a consequence, it may not be Hausdorff, depending on the properties of $\phi$. We avoid this problem as follows. Choose $c_0$ so that $c_U$ is a point. Instead of using in Theorem \ref{foliation} a subspace $L$ of minimal codimension as we did, we take $L'$ such that
$$
L=L'\oplus L_1
$$
where $L_1$ is one-dimensional. Then using Theorem \ref{foliation} with $L'$ this time, we obtain in place of $K_U$ the space $K_U\times W_1$ for some open set $W_1\subset L_1$. Now, we may assume that the image of $c_0$ in $K_U\times W_1$ is not constant (for example that the projection onto $L_1$ is not constant). Taking this image as the new path $c_U$, we may now finish the argument with this $c_U$, defining $K_{c_0}$ as the gluing of $K_U\times W_1$ through $\phi$ as before. Observe that the gluing is given by the associated map
$$
(J,v)\in K_U\times W_1\longmapsto (\phi\cdot J,d\phi\cdot v)\in K_U\times W_1
$$
the action on the first coordinate being defined in \eqref{action}, and the action on the second coordinate being that of the differential of $\phi$ on $\Sigma(TX^{diff})$.
\vspace{5pt}\\
Notice that, when dealing with $L$ of minimal codimension, that is excluding the case where $\mathcal X_0$ is a non-trivial CR bundle, it follows from Theorem \ref{foliation} and Theorem \ref{foliationtrivial}
that $K_U$, hence also $K_{c_0}$, is complete at $c_U(t)$ but not always versal, being the product of the Kuranishi space with a $\mathbb C$-vector space of dimension 
$$
\text{codim } L-h^0(t).
$$
Hence it has dimension
$$
h^1(t)+\text{codim } L- h^0(t)
$$
as stated. In the case of a non-trivial CR bundle, we have to increase the dimension by one. 
\vspace{5pt}\\
To define the map $\Xi_\phi$, we proceed as follows. Assume that we are not in the special case. We already have an a.c. projection
\begin{equation}
\label{projection}
\begin{CD}
U\cap \mathcal I@ > a.c. >> K_U
\end{CD}
\end{equation}
by Theorem \ref{foliationtrivial}. Since both the construction of $U_\phi$ and that of $K_{c_0}$ consist in taking the quotient by $\phi$, it follows that the projection \eqref{projection} descends as a map
$$
\begin{CD}
U_\phi\cap\mathcal I@ >\Xi_\phi>> K_{c_0}
\end{CD}
$$
as desired and that $K_{c_0}$ is the leaf space of the $L$-foliation of $U_\phi$. The special case is handled in the same way, just noting that, running the proof of Theorem \ref{foliationtrivial} with $L'$ instead of $L$ yields an a.c. projection
\begin{equation}
\label{projectionbis}
\begin{CD}
U\cap \mathcal I@ > a.c. >> K_U\times W_1.
\end{CD}
\end{equation}
Using \eqref{projectionbis} instead of \eqref{projection}, we immediately see that it descends also as a map \eqref{pathretraction}.
\end{proof}

\begin{remark}
If $K_U$  is not reduced, then we perform exactly the same construction. Because of \eqref{tildeLi}, the extra vector space factor is always reduced and there is no change in the dimension counting.
\end{remark}

\begin{corollary}
\label{diffcorollary}
The map {\rm \eqref{pathretraction}} is, after shrinking of $U_\phi$, a.c. diffeomorphic to a trivial bundle with base $K_{c_0}$ and fiber a submanifold of $\text{\rm Diff}^0(X^{diff})$ passing through the identity.
\end{corollary}

\begin{proof}
Use trivialization \eqref{localcharts} and observe that the $\phi$-gluing respects the fibers of this trivialization. This shows that $U_\phi$ is a.c. diffeomorphic to a locally trivial bundle over $K_{c_0}$. Now there are two cases. Either the gluing occurs on the fibers (cf. case $1$ and $2$ in the proof of Theorem \ref{foliationtrivial}) hence the bundle is trivial; or the gluing occurs on the base (cf. case $3$ of the proof of Theorem \ref{foliationtrivial}), but then its monodromy is isotopic to the identity (because it is given by some map obtained as $e(\eta)$ once the gluing between $U_k$ and $U_1$ is performed), and it is trivial.
\vspace{5pt}\\
Observe that, in the first case, the fiber is homotopic to a circle, whereas in the second case it is contractible but this time the base has a non-trivial fundamental group.
\end{proof}

We also note the
\begin{corollary}
\label{corollary4}
The analytic space $K_{c_0}$ is unique in the following sense.
\begin{enumerate}
\item[(i)]  Up to a.c. isomorphism, it does not depend on the choice of the adapted covering and of the Kuranishi maps.
\item[(ii)] If $U'$ is another neighborhood of $c_0$, then the restrictions of $\Xi_\phi$ and of the corresponding map $\Xi'_\phi$ to $(U\cap U')_\phi\cap \mathcal I$ have a.c. isomorphic images.
\end{enumerate}
\end{corollary}

\begin{proof}
This follows immediately from the fact that it is a leaf space.
\end{proof}

In other words, the germ of $K_{c_0}$ at $c_0$ is unique. However, once again, we want to emphasize that it depends on the choice of $L$, and of $L'$ in the special case.

\section{Deformations in the spirit of Kodaira-Spencer.}
\label{deformations}

 \subsection{Levi-flat CR spaces}
 \label{LeviCRss}
 Recall definition \ref{deformationDef}. Note that it holds for infinite-dimensional (Banach or, by extension and abuse of notation, Fr\'echet) analytic spaces as bases. For example, $\mathcal I$ is the base of such an infinite-dimensional analytic deformation of $X_0$, once a base-point corresponding to $X_0$ is fixed. This comes from the fact that, using the map \eqref{Kurretraction} at each point of $\mathcal I$, one shows that this family is locally obtained by pull-back from the Kuranishi family. Since the last one is flat, so is the first one.\vspace{5pt}\\
We need a notion of flat family for proper CR submersions onto the circle. The basic idea is the following one. Let $\mathcal X_0$ be a proper CR submersion and let $X_0$ be a fiber of $\mathcal X_0$. If $\mathcal Z\to B$ is a deformation of  $\mathcal X_0$, then $\mathcal Z$ should admit a submersion onto $B\times\mathbb S^1$ making it a deformation of the compact complex manifold $X_0$ with base $B\times\mathbb S^1$. However, this induces some technical difficulties, because $\mathcal Z$ is neither a smooth manifold nor an analytic space. When $B$ is an analytic space, it is not even a CR manifold but a sort of singular CR object that we define now. The corresponding notion of flatness, that of a transflat morphism is then defined in \ref{transflatss}. Finally, we give the definition in subsection \ref{DefCRss}.

 \begin{definition}
A {\it Levi-flat CR space} $Z$ is a second-countable Hausdorff space for which there exists a covering by open subsets $V_\alpha$ and homeomorphisms
$$
F_\alpha \ :\ V_\alpha\longrightarrow \mathbb R^p\times W_\alpha
$$
for some analytic sets $W_\alpha\subset\mathbb C^{n_\alpha}$, such that the changes of charts
$$
F_{\alpha\beta}:=F_\beta\circ F_\alpha ^{-1}
$$
are smooth, respect the foliation by copies of $W_\alpha$, and are analytic in the second variable; that is, we have
$$
F_{\alpha\beta}\ :\ (x,z)\longmapsto (f_{\alpha\beta}(x),g_{\alpha\beta}(x,z))
$$
and, for all $x$, the map
$$
z\longmapsto g_{\alpha\beta}(x,z)
$$
is analytic.
\end{definition}

In this definition, the analytic sets $W_\alpha$ may be non-reduced. Alternatively, we could have described a Levi-flat CR space as a ringed space locally homeomorphic to $\mathbb R^p\times W_\alpha$ with structure sheaf given by the sheaf of functions analytic on $W_\alpha$ for all $x\in\mathbb R^p$ and smooth in $x$. A Levi-flat CR space belongs to the special class of ringed spaces called mFB spaces (see \cite[p. 119]{F-K}). But we stick to a more geometric definition, since we really want to consider it as a Levi-flat CR manifold with singularities. As in the smooth case, it is foliated, the leaves being obtained by gluing the $W_\alpha$ via $g_{\alpha\beta}$. The leaves are analytic, but, unlike the smooth case, they may have singularities. A trivial example is given by a product of an analytic space with a smooth manifold. Especially, $B\times\mathbb S^1$ is a Levi-flat CR space for any analytic space $B$.

\subsection{Transflat morphisms}
\label{transflatss}
We define
\begin{definition}
A CR morphism $\Pi : \mathcal Z\to B$ between Levi-flat CR spaces is {\it transflat} if there are submersion charts 
$$
\begin{CD}
z\in U\subset \mathcal Z @ >\text{CR iso.}>> \Pi (U)\times \mathbb C^n\times\mathbb R^p\cr
@ V \Pi VV @ VV\text{1st projection}V\cr
\Pi(z)\in\Pi (U)\subset B @ >>Id>   \Pi(z)\in\Pi (U)\subset B
\end{CD}
$$
 for all points $z\in\mathcal Z$ (see \cite{Sc}).
 \end{definition}
 
\begin{remark}
In the previous definition, if the fibers of $\mathcal Z$ are complex manifolds, we have $p=0$ in the diagram of submersion charts. If $B$ is an anlytic space, we get exactly the geometric definition of a flat family (with smooth fibers).
\end{remark}

From this definition, we have
\begin{definition}
A {\it CR deformation} of $X_0$ is a Levi-flat CR space $\mathcal Z$ together with a proper and transflat CR morphism $\Pi : \mathcal Z\to B$, for $B$ a Levi-flat CR space, a base-point and a marking.
\end{definition}

In \cite{F-K} and \cite{Sc}, such a deformation is called a "relativ-analytisch Deformazion". Observe that analytic and smooth deformations are particular cases of CR deforma\-tions. Observe also that if $B$ is a product $B_1\times B_2$ with $B_1$ smooth and $B_2$ analytic, then for every $x\in B_1$, the induced deformation over $B_2$ obtained by restricting $\mathcal Z$ to $\Pi ^{-1}(\{x\}\times B_2)$ is analytic; whereas for every $z\in B_2$, the induced deformation over $B_1$ obtained by restricting $\mathcal Z$ to $\Pi ^{-1}(B_1\times\{z\})$ is smooth.\vspace{5pt}\\
Observe that, if $f : B\to \mathcal I$ is CR, then the family $(J_{f(t)})_{t\in B}$ defines a CR deformation of $X$.

\subsection{Deformations of proper CR submersions}
\label{DefCRss}

Let $\pi : \mathcal X_0\to \mathbb S^1$ be a proper CR submersion. The following definition is inspired in \cite{Bu}.

\begin{definition}
\label{DefCRDef}
A {\it holomorphic deformation} (resp. {\it smooth deformation}) of $\mathcal X_0$ is a Levi-flat CR space $\mathcal Z$ together with a proper and transflat CR morphism $\Pi : \mathcal Z\to B$ onto an analytic space (resp. a smooth manifold) $B$, a smooth and proper transflat map $s : \mathcal Z\to \mathbb S^1$ and a marking $i : \mathcal X_0\to \Pi^{-1}\{0\}$ such that, for all $b\in B$,
\begin{enumerate}
\item[(i)] The $\Pi$-fiber $\mathcal Z_b$ over $b$ is a Levi-flat CR submanifold of $\mathcal Z$.
\item[(ii)] The restriction $s_b$ of $s$ to $\mathcal Z_b$ is a proper CR submersion compatible with $\pi^{diff}$.
\item[(iii)] The composition $s_0\circ i$ is equal to $\pi$.
\item[(iv)] The map
\begin{equation}
\label{Pmap}
P=(s,\Pi)\ :\ \mathcal Z\to \mathbb S^1\times B
\end{equation}
is a proper and transflat CR morphism.
\end{enumerate}
\end{definition}

This is a quite technical definition so let us highlight some of its principal features. Firstly, choose any $t_0\in\mathbb S^1$ and let $X_0$ be the fiber of $\mathcal X_0$ over $t_0$. Then the map $P$ of \eqref{Pmap}
is a CR deformation of $X_0$ once chosen a marking.\vspace{5pt}\\
Secondly, $\mathcal Z$ is locally diffeomorphic at $z\in\mathcal Z$ to $\mathcal X^{diff}\times U$, for $U$ a neighborhood of $\Pi (z)$ in $B$. Moreover, we have a commutative diagram of {\it diffeomorphisms}
$$
\begin{CD}
\mathbb S^1@ <\pi^{diff}<< \mathcal X^{diff}\cr
@ A s AA @ AA\text{1st projection}A\cr
\Pi^{-1}(U)\subset\mathcal Z@ >\simeq >>\mathcal X^{diff}\times U \cr
@ V\Pi VV  @ VV \text{2nd projection} V\cr
U @ >>\text{Identity}> U\cr
\end{CD}
$$
Thirdly, since $\Pi$ and $P$ are transflat, in the neighborhood of any point of $\mathcal Z$, we have CR coordinates $(z,t,b)$ with $z$ holomorphic coordinates of the compact complex manifold $P^{-1}\{(t,b)\}$, and $b$ of the base $B$, whereas $t$ is given by the value of $s$.\vspace{5pt}\\
Fourthly, choose any $t\in\mathbb S^1$. Define
\begin{equation}
\label{Yt}
\begin{CD}
\Pi_t\ :\ \mathcal Y_t:=s^{-1}(t)@> \restriction{\Pi}{\mathcal Y_t}>> B.
\end{CD}
\end{equation}
This is a deformation of the compact complex manifold $X_t:=\pi^{-1}(t)$.\vspace{5pt}\\
We may easily define isomorphism of deformations of proper CR submersions (which we assume to induce the identity on the base) and pull-back. We omit the details.

\section{Kodaira-Spencer map.}
\label{KSmapsection}

Let us start with the construction of the Kodaira-Spencer map of a deformation of a proper CR submersion over the circle. In the classical case, the Kodaira-Spencer map of a deformation of $X_0$ takes value in the first cohomology group $H^1(X_0,\Theta_0)$, which can be identified with the tangent space at $0$ of the Kuranishi space of $X_0$ in such a way that it corresponds to the differential $d_0f$ of the map $f$ obtained by completeness. In our case, however, the Kodaira-Spencer map will take value in a first cohomology group which is different from $T_0K^g$, see section \ref{versalitysection}.\vspace{5pt}\\
To do that, we start as usual with a proper CR submersion $\pi : \mathcal X_0\to\mathbb S^1$ and we define $\Theta_\pi$ to be the sheaf of germs of CR vector fields of $\mathcal X_0$ tangent to the fibers of $\pi$. The first cohomology group $H^1(\mathcal X_0,\Theta_\pi)$ has a natural structure of the set of smooth sections of a sheaf of $\mathbb C$-vector spaces over the circle. The stalk at some point $t$ is the vector space $H^1(X_t,\Theta_t)$. Moreover, it is a vector bundle over the circle as soon as the function $h^1$ is constant along the circle, cf. \cite{K-S1}.\vspace{5pt}\\
Let $\Pi : \mathcal Z\to B$ be a deformation of $\mathcal X_0$. Let $\Theta_s$ be the sheaf of germs of CR vector fields of $\mathcal Z$ tangent to the fibers of $s$ and whose image through the differential of $\Pi$ is a germ of vector field of $B$. In CR coordinates $(z,t,b)$ (see the end of section \ref{KSmapsection}), such a germ is written as
\begin{equation}
\label{localform}
\sum a_i(z,t,b)\dfrac{\partial}{\partial z_i}+\sum c_j(b)\dfrac{\partial}{\partial b_j}.
\end{equation}
Let also $\Theta_P$ be the sheaf of germs of CR vector fields of $\mathcal Z$ tangent to the fibers of $P$. Consider the following exact sequence of sheaves
\begin{equation}
\label{fundseq}
\begin{CD}
0@ >>> \Theta_P @ >>>\Theta_s@ >>> \Theta_s/\Theta_P@ >>>0
\end{CD}
\end{equation}
and observe that the quotient sheaf $\Theta_s/\Theta_P$ can be identified with the sheaf $\Theta_B$ of germs of CR vector fields on the base $B$ because of \eqref{localform}. The long exact sequence associated to \eqref{fundseq} runs as follows
\begin{equation}
\label{longfundseq}
\begin{CD}
\hdots\  H^0(\mathcal Z,\Theta_s) @ >>> H^0(B,\Theta_B)@ >\rho >>H^1(\mathcal Z,\Theta_P)@ >>>\hdots
\end{CD}
\end{equation}
Observe now that the restriction of $\rho$ to the tangent space $T_0B$ gives a map
\begin{equation}
\label{KSmap}
\begin{CD}
T_0B@ >\rho_0>> H^1(\mathcal X_0,\Theta_\pi).
\end{CD}
\end{equation}
\begin{definition}
The map $\rho_0$
of {\rm \eqref{KSmap}} is called the {\it Kodaira-Spencer map} at $0$ of the deformation $\Pi : \mathcal Z\to B$.
\end{definition}
Roughly speaking, the map $\rho$ of \eqref{longfundseq} represents the complete obstruction to lift CR vector fields of $B$ to CR vector fields of $\mathcal Z$ respecting the fibers of $s$, thus trivializing the family. Indeed, \eqref{longfundseq} is the exact analogue of the fundamental exact sequence of \cite{K-S1}. The Kodaira-Spencer map being the evaluation of this obtruction to the central fiber is the first obtruction to such a trivialization. It has the advantage to be defined on $\mathcal X_0$ and not on the whole deformation $\mathcal Z$ and thus can be computed explicitely in many cases.

\section{Completeness.}
\label{completesection}
We are in position to prove the following {\it completeness} result.

\begin{theorem} 
\label{completetheorem}
Let $\pi : \mathcal X_0\to\mathbb S^1$ be a proper CR submersion. Choose a path $c_0$ representing $\mathcal X_0$ and satisfying Hypothesis {\rm \ref{noselfint}}. Let $K^g$ be a Kuranishi type moduli space of $\mathcal X_0$.\smallskip

Then, there exists a (infinite-dimensional) holomorphic deformation $\Pi^g : {\mathcal K}^g\to {K}^g$ of $\mathcal X_0$ with base ${K}^g$ which is complete at $0$. Moreover, we may ask the local isomorphisms given by completeness to preserve the markings.
\end{theorem}

\begin{proof}
Choose a path $c_0$ representing $\mathcal X_0$, define $K_{c_0}$ and $K^g$ as usual. 
Observe $K_{c_0}$ being obtained from a finite number of Kuranishi spaces by gluing them through a.c. isomorphisms, it defines in a natural way an analytic family $\mathcal K_{c_0}$.\vspace{5pt}\\
Indeed, going back to the proof of Theorem \ref{foliationtrivial} and using the notations introduced there, let $K_1$ and $K_2$ be two Kuranishi spaces and let 
\begin{equation}
\label{K12}
K_{12}=K_1\cup_g K_2
\end{equation}
be an analytic space obtained by gluing, the gluing $g$ being an a.c. isomorphism between an open set $V_1$ of $K_1$ and an open set $V_2$ of $K_2$. More precisely, $g$ is given as a map
$$
g\ :\ J\in V_1\subset K_1\longmapsto G(J)\cdot J\in K_2
$$
where 
$$
G\ :\ V_1\longrightarrow \text{Diff}^0(X^{diff})
$$
and where the $\cdot$ denotes the action \eqref{action} (this $G$ is nothing else than the $e(\eta)$ of \eqref{geta}). Consider now the families $\mathcal K_i$ induced above $K_i$. They are constructed as $K_i\times X^{diff}$, every fiber $\{J\}\times X^{diff}$ being endowed with the complex structure $J$. It follows that the map
\begin{equation}
\label{gluingg}
(J,x)\in V_1\times X^{diff}\longmapsto (g(J),G(J)(x))\in V_2\times X^{diff}
\end{equation}
realizes an analytic isomorphism between open sets of $\mathcal K_i$ preserving the projections onto $K_i$. Gluing $\mathcal K_1$ and $\mathcal K_2$ through \eqref{gluingg}, we obtain a holomorphic family $\mathcal K_{12}$ over $K_{12}$. Repeating the process yields the family $\mathcal K_{c_0}$.\vspace{5pt}\\
Set 
$$
{\mathcal K}^g:=\{f^*\mathcal K_{c_0}\mid f\in C^\infty(\mathbb S^1, K_{c_0})\}
$$
and let $\Pi^g$ be the induced projection from $\mathcal K^g$ to $K^g$. Call $0$ the point of $K^g$ corresponding to $c_0$. Finally choose a marking 
$$
i^g \ :\ \mathcal X_0 \longrightarrow (\Pi^g)^{-1}(0).
$$
This is obviously a holomorphic deformation of $\mathcal X_0$. Moreover,
as a consequence of property \eqref{Pmap}, a holomorphic (resp. smooth) deforma\-tion $\Pi : \mathcal Z\to B$ of $\mathcal X_0$ can be locally encoded as follows. Construct the neighborhood $U_\phi$ of $c_0$. Then $\Pi$ induces an analytic (resp. smooth) map $f$ from a neighborhood of $0$ in $B$ to $C^\infty (\mathbb S^1, U_\phi\cap\mathcal I)$ such that $\mathcal Z$ is locally isomorphic to $\mathcal X^{diff}\times B$ endowed with the family of CR submersions structures $(f(t))_{t\in B}$.\vspace{5pt}\\
It is enough to compose $f$ with the a.c. projection from $C^\infty (\mathbb S^1, U_\phi\cap\mathcal I)$ to $K^g$ induced by the map  $\Xi_\phi$ of \eqref{pathretraction} to prove completeness.\vspace{5pt}\\
Finally, observe that if this composition, let us call it $g$, does not preserve the markings, then there exists a CR isomorphism $\Phi$ of $\mathcal X_0$ covering the identity of $\mathbb S^1$ such that
$$
\begin{CD}
\mathcal X_0@ >\Phi>>\mathcal X_0\cr
@ ViVV @ VVi^gV\cr
(\mathcal Z,\Pi ^{-1} (0)) @ >G>>\mathcal K^g\cr
@ V\Pi VV @ VV \Pi^g V\cr
(B,0) @ >>g >K^g
\end{CD}
$$
But now, $\Phi$ being an element of $\text{Diff}^0_\pi(\mathcal X^{diff})$ acts on $\mathcal K^g$, so that we may replace the map $G$ by $\Phi^{-1}\cdot G$ so that we have now
$$
\begin{CD}
\mathcal X_0@ > i^g >> \mathcal K^g\cr
@ ViVV @ VV\text{Identity} V\cr
(\mathcal Z,\Pi^{-1}(0))@ >>\Phi^{-1}\cdot G>\mathcal K^g
\end{CD}
$$
and this time the map $g$ preserves the markings.
\end{proof}
\begin{remark}
Although the map 
$$
\Pi^g\ :\ \mathcal K^g\longrightarrow K^g
$$
is not strictly speaking a deformation of $\mathcal X_0$ (since its base is only Fr\'echet), we will consider it as a deformation of $\mathcal X_0$. The reader may, if he wants, replace it by its Sobolev completion in the sequel.
\end{remark}

\begin{remark}
\label{pullbackdefremark}
Observe that, if $B$ is a smooth manifold, respectively an analytic space, with base point $0$ and if 
$$
f\ :\ B\times\mathbb S^1\longrightarrow K_{c_0}
$$
is a smooth map, respectively a CR morphism, sending $\{0\}\times\mathbb S^1$ onto $c_0$, then $f^*\mathcal K_{c_0}$ induces naturally a deformation of $\mathcal X_0$ once chosen a marking (compare with subsection \ref{EandIss}). This a simple way to construct deformations of proper CR submersions and we will often use this trick in the sequel.
\end{remark}

\section{(Uni)versality.}
\label{versalitysection}
As recalled in subsection \ref{classicalKSss}, the Kuranishi family is not only complete but also versal at $0$ in the classical case. We will see that things are a bit different in the case of proper CR submersions.\vspace{5pt}\\
First, recall that we gave two equivalent definitions of versality in the classical context. The first one deals with the Kuranishi space having minimal dimension at $0$. Since our family is infinite-dimensional, we cannot use this definition. The other characterization deals with the uniqueness of the differential at $0$ of the map $f$ associated by completeness to a deformation. This can be easily transposed to our context.\vspace{5pt}\\
Moreover, recall that the Kuranishi space is {\it universal} if the germ of $f$ is unique. This definition can also be transposed to our context. It is known that, in the reduced case, the Kuranishi space is universal if and only if the function $h^0$ is constant along it (cf. \cite{Wa1}, \cite{Wa2}, \cite{Me1}).\vspace{5pt}\\ 
Let us fix some notations. We start from a deformation $\Pi : \mathcal Z\to B$ of $\mathcal X_0$ and we consider a map $f$ from some open set of $B$ to $K^g$ given by completeness. Its differential at the marked point $0$ goes from the tangent space $T_0B$ to the tangent space $T_0K^g$. We only consider local isomorphisms from $\mathcal Z$ to $\mathcal K^g$ covering $f$ which preserve the markings. This prevents from composing with automorphisms of $\mathcal Z$ and $\mathcal K^g$ which descend as non-trivial automorphisms of the bases of the families. Notice that this is necessary to hope the uniqueness of $d_0f$. Even in the classical case, allowing these compositions, one loses versality.\vspace{5pt}\\
Nevertheless, we do not have versality in general. We will give some counterexamples in section \ref{exemples} and Remark \ref{versalremark}. Indeed,

\begin{theorem}
\label{versaltheorem}
We have
\begin{enumerate}
\item[(i)] The family $\Pi^g : \mathcal K^g\to K^g$ is versal at $0$ if and only if
\begin{equation}
\label{versalcondition}
\forall t\in\mathbb S^1,\qquad h^0(t)=\text{\rm codim } L.
\end{equation}
\item[(ii)] The family $\Pi^g : \mathcal K^g\to K^g$ is universal at $0$ for families over a reduced base if and only if
\begin{equation}
\label{universalcondition}
\forall J\in K_{c_0},\qquad h^0(X_J)=\text{\rm codim } L.
\end{equation}
\end{enumerate}
\end{theorem}

\begin{remark}
\label{versalremark}
Observe that, if $\mathcal X_0$ is a trivial CR bundle, then \eqref{versalcondition} is automatically satisfied, hence $\mathcal K^g$ is versal. And it is universal if and only if the Kuranishi space of the fiber is universal. On the contrary, if $\mathcal X_0$ is a {\it non-trivial} CR bundle, the family $\mathcal K^g$ is never versal.
\end{remark}

\begin{proof}
We first need to compute the tangent space of $K^g$ at some point $c$. From its definition, one has that
$$
T_cK^g=\{H\in C^\infty (\mathbb S^1, TK_{c_0})\mid p\circ H\equiv c\}
$$
where $p : TK_{c_0}\to K_{c_0}$ is the tangent sheaf of $K_{c_0}$. In particular, for $t\in\mathbb S^1$, we have a commutative diagram
\begin{equation}
\label{evaCD}
\begin{CD}
H\in T_cK^g@ >\text{evaluation at }t>>H(t)\in T_{c(t)}K_{c_0}\cr
@ V p\circ VV @ VV p V\cr
c\in K^g@ >\text{evaluation at }t>>(p\circ H)(t)=c(t)\in K_{c_0}
\end{CD}
\end{equation}
Let $ev_t$ denote the evaluation map of the bottom arrow and $Ev_t$ that of the top arrow. Then analyzing \eqref{evaCD} yields that the differential $d_cev_t$ is equal to $Ev_t$.\vspace{5pt}\\
Indeed, let $v$ be a vector of $T_cK^g$ and let
$$
u\ :\ (-\epsilon,\epsilon)\longrightarrow K^g
$$
be a smooth path whose derivative at $0$ is $v$. Define
$$
(s,t)\in (-\epsilon,\epsilon)\times \mathbb S^1\longmapsto U(s,t):=ev_t \circ u(s)\in K_{c_0}.
$$
Compute
$$
\restriction{\dfrac{d}{ds}}{s=0} (ev_t\circ u(s))=\restriction{\dfrac{d}{ds}}{s=0}(U(s,t))=Ev_t\circ \restriction{\dfrac{d}{ds}}{s=0}(u(s))
$$
that is
$$
d_c(ev_t\circ u)=d_c ev_t (v)=Ev_t\circ v.
$$

Let $\Pi : \mathcal Z\to B$ be a deformation of $\mathcal X_0$ and let $(f,F)$ be given by completeness, i.e. $F$ is the isomorphism between $\mathcal Z$ and $f^*\mathcal K^g$. Consider now the deformation of $X_t:=\pi^{-1}(t)$ defined in \eqref{Yt}
$$
\begin{CD}
\Pi_t\ :\ \mathcal Y_t:=s^{-1}(t)@> \restriction{\Pi}{\mathcal Y_t}>> B.
\end{CD}
$$
We have a commutative diagram
\begin{equation}
\label{CDbis}
\begin{CD}
\mathcal Y_t\subset\mathcal Z@ >(F)_t>>\mathcal K_{c_0}\cr
@ V \Pi_t VV @ VVV\cr
(B,0)@ >> (f)_t> (K_{c_0}, c_0(t))
\end{CD}
\end{equation}
where $(F)_t$ and $(f)_t$ are given by the evaluation at $t$ of the maps $F$ and $f$.\vspace{5pt}\\
Observe that, if
$$
i\ :\ \mathcal X_0\longrightarrow \Pi^{-1}(0)
$$
denotes the marking of our family $\Pi$, then 
$$
i_t\ :\ X_t\subset \mathcal X_0\longrightarrow \Pi_t^{-1}(0)
$$
is a marking for $\Pi_t$.\vspace{5pt}\\
 From \eqref{evaCD} and \eqref{CDbis}, we obtain that
\begin{equation}
\label{dfequal}
d_0(f)_t\equiv (d_0f)_t\ :\ T_0B\longrightarrow T_{c_0(t)}K_{c_0}.
\end{equation}
and we see that the versality at $c_0(t)$ of $K_{c_0}$ for any $t$ implies the versality at $0$ of $K^g$. The markings used here are $i_t$ and $i$.\vspace{5pt}\\
Conversely, assume that $K_{c_0}$ is not versal at some point $c_0(t_0)$ with $t_0\in\mathbb S^1$. Then we can find a deformation
$$
\mathcal Y\to (\mathbb D,0)
$$
of $X_{t_0}$ with marking $i_0$ and two holomorphic maps
$$
\begin{CD}
(\mathbb D,0)@ >f_{t_0},\ g_{t_0}>> (K_{c_0},c_0(t_0))
\end{CD}
$$
with
$$
\mathcal Y=f_{t_0}^*\mathcal K_{c_0}=g_{t_0}^*\mathcal K_{c_0},
$$
respecting the markings $i_0$ and $i^g(c_0(t_0))$ and finally such that 
\begin{equation}
\label{notequiv}
d_0f_{t_0}\not\equiv d_0g_{t_0}.
\end{equation}
The proof will consist in extending $f_{t_0}$ and $g_{t_0}$ into CR maps $f$ and $g$ over $\mathbb D\times\mathbb S^1$ in such a way that $f^*\mathcal K_{c_0}$ and $g^*\mathcal K_{c_0}$ will define exactly the same deformation of $\mathcal X_0$ but with $f$ different from $g$.\vspace{5pt}\\
Firstly, we deal with the extension of the map $f_{t_0}$ into a CR map
$$
f\ :\ \mathbb D\times\mathbb S^1\longrightarrow K_{c_0}
$$
such that
$$
\left\{
\begin{aligned}
&\restriction{f}{\mathbb D\times\{t_0\}}\equiv f_{t_0}\cr
&\restriction{f}{\{0\}\times\mathbb S^1}\equiv c_0
\end{aligned}
\right .
$$
To prove that such an extension exists, first observe that we may assume $K_{c_0}$ smooth. If not, just desingularize and lift both $f_{t_0}$ and $c_0$. Also observe that we may perform the extension step by step along the circle, from $t_0$ to a close $t_1$ and so on. Hence we reduce the problem to a local extension problem in $\mathbb C^n$.

\begin{remark}
There is no particular problem if the space $K_{c_0}$ is not reduced. Indeed, since we use maps from a reduced base (here a circle or an annulus) into $K_{c_0}$, they map into the reduction of $K_{c_0}$. Hence, we only have to desingularize the reduction of $K_{c_0}$. 
\end{remark}

Taking into account that we have a monodromy problem when coming back to $t_0$ after a complete turn, we finally see that the proof is completed with the following lemma (we prove more since it will be useful in section \ref{rigiditysection}).

\begin{lemma}
\label{extensionlemma}
We have
\begin{enumerate}
\item[(i)] Let $u : [0,1]\to\mathbb C^n$ be a smooth path and let $f_i : \mathbb D\to\mathbb C^n$ be two ($i=0,1$) holomorphic maps such that $f_i(0)=u(i)$.
Then, there exists a CR map
$$
F\ :\ \mathbb D\times [0,1]\longrightarrow \mathbb C^n
$$
such that 
$$
\left\{
\begin{aligned}
F(-,i)&\equiv f_i\quad\text{ for }i=0,1\cr
F(0,t)&=u(t)\quad\text{ for }t\in [0,1]
\end{aligned}
\right .
$$
\item[(ii)] Moreover, let $v : [0,1]\to \mathbb C^n$ be a smooth path and $z_0\in\mathbb D\setminus\{0\}$ such that $f_i(z_0)=v(i)$ (for $i=0,1$).
Then, we may assume that the map $F$ of (i) satisfies also
$$
F(z_0,t)=v(t)\quad\text{ for }t\in [0,1].
$$
\end{enumerate}
\end{lemma}

\begin{proof}[Proof of Lemma 2]
To prove (i), just define
$$
F(z,t)=(1-t)f_0(z)+tf_1(z)+u(t)-tu(1)-(1-t)u(0),
$$
and to prove (ii),
$$
\begin{aligned}
F(z,t)=&(1-t)f_0(z)+tf_1(z)+u(t)-tu(1)-(1-t)u(0)\cr
+&\dfrac{z}{z_0}(v(t)-(1-t)v(0)-tv(1)-u(t)+tu(1)+(1-t)u(0)).
\end{aligned}
$$
\end{proof} 
We are in position to define
\begin{equation}
\label{Zdef}
\begin{CD}
\mathcal Z:=f^*(\mathcal K_{c_0})\longrightarrow \mathbb D\times\mathbb S^1@ >\text{1st projection} >>\mathbb D.
\end{CD}
\end{equation}
This is a deformation of $\mathcal X_0$ with marking 
$$
i\ :\ \mathcal X_0\longrightarrow c_0^*\mathcal K_{c_0}.
$$
Secondly, we deal with the extension of $g_{t_0}$ into a CR map from $\mathbb D\times\mathbb S^1$ such that $g^*\mathcal K_{c_0}$ also defines $\mathcal Z$. To do this, we proceed in a completely different way. We need to reinterpret the construction of $\mathcal Y$ and $\mathcal Z$.
Taking into account that $K_{c_0}$ encodes complex operators, we may rewrite $\mathcal Y$ as
$$
\mathcal Y=(X^{diff}\times \mathbb D, J_0)=(X^{diff}\times \mathbb D, J_1)
$$
with $J_0$ and $J_1$ families of complex operators on $X^{diff}$ indexed by $\mathbb D$ and satisfying
$$
(F_{t_0})_*J_0\equiv (G_{t_0})_*J_1
$$
where $F_{t_0}$ (respectively $G_{t_0}$) is a map from $\mathcal Y$ to $\mathcal K_{c_0}$ covering $f_{t_0}$ (respectively $g_{t_0})$.\vspace{5pt}\\
In other words, we may find some smooth family of diffeomorphisms $k_s$ of $X^{diff}$ parametrized over the disk such that, for every $s\in\mathbb D$, the map $k_s$ is a biholomorphism from $(X^{diff}\times\{s\}, J_0(s))$ to $(X^{diff}\times\{s\}, J_1(s))$.\vspace{5pt}\\
The previous extension of $f_{t_0}$ into $f$ is nothing else than extending $J_0$ into a family of complex operators over $\mathbb D\times \mathbb S^1$ so that 
$$
\mathcal Z=(X^{diff}\times \mathbb D\times\mathbb S^1, J_0)
$$
To define $g$, we extend $J_1$ over the same base by defining
$$
J_1:=(k_s)_*J_0
$$
on the fiber over any point $(s,t)\in\mathbb D\times\mathbb S^1$.\vspace{5pt}\\
Obviously, we also have
$$
\mathcal Z=(X^{diff}\times \mathbb D\times\mathbb S^1, J_1).
$$
Shrinking the base $\mathbb D$ if necessary, we may assume that the map
$$
J_1\ :\ \mathbb D\times\mathbb S^1\longrightarrow \mathcal I
$$
has image in the open set $U_\phi\cap\mathcal I$ admitting a retraction \eqref{pathretraction} onto $K_{c_0}$. This allows us to extend the map $g_{t_0}$ into a map $g$ defined over $\mathbb D\times\mathbb S^1$ by stating
$$
g(z,t):=\Xi_\phi ((J_1){(z,t)})
$$
where $(J_1){(z,t)}$ denotes the restriction of $J_1$ to the fiber over $(z,t)$. By construction, we have
$$
\mathcal Z=g^*\mathcal K^g.
$$
Because of \eqref{dfequal}, \eqref{notequiv} and \eqref{Zdef}, the family $\mathcal K^g$ is not versal at $c_0$.\vspace{5pt}\\
To finish with the proof of (i), we just have to show that $K_{c_0}$ is versal at each point $c_0(t)$ if and only if equality \eqref{versalcondition} holds.\vspace{5pt}\\
The construction of $K_{c_0}$ given in the proof of \ref{foliationtrivial} shows that it is complete at $c_0(t)$ with dimension 
\begin{equation}
\label{dimen}
h^1(t)+\text{codim } L-h^0(t)\quad\text{ or }\quad h^1(t)+\text{codim } L+1-h^0(t).
\end{equation}
To be versal, it must have minimal dimension, that is dimension $h^1(t)$. Since we have
$$
\text{codim } L\geq h^0(t) \qquad\text{for all }t
$$
this yields the condition \eqref{versalcondition}.\vspace{5pt}\\
Conversely, if this condition is fulfilled, then, from \eqref{dimen}, the space $K_{c_0}$ is versal at each point of $c_0$, hence the differential $d_0(f)_t$ is uniquely determined for each $t\in\mathbb S^1$. By \ref{dfequal}, which means that $d_0f$ is uniquely determined, so the family $\Pi^g : \mathcal K^g\to K^g$ is versal at $0$. This proves (i).\vspace{5pt}\\
By definition $f$ is unique, yielding universality if $f_t$ is unique for all $t$. In other words, $K^g$ is universal if  $K_{c_0}$ is universal at each point $c_0(t)$. Conversely, if $K_{c_0}$ is not universal at some point $c_0(t_0)$, then the same argument as above (just replacing \eqref{notequiv} with $f_{t_0}\not\equiv g_{t_0}$) shows that $K^g$ is not universal at $c_0$.\vspace{5pt}\\
Universality of $K_{c_0}$ implies also that it is versal at each point $c_0(t)$, that is that the function $h^0$ is equal to $\text{codim } L$ at every point $c_0(t)$. Now, by a Theorem of Wavrik (see \cite[\S 5.5]{Me1} for the version we use), $K_{c_0}$ is universal at $c_0(t)$ for families over a reduced base if and only if $h^0$ is constant in a whole neighborhood of $c_0(t)$. So we finally obtain that the condition
$$
h^0(X_J)=\text{codim } L\qquad\text{ for all }J\in K_{c_0}
$$
is sufficient to have universality for families over a reduced base.
\end{proof}

\section{Back to the Kodaira-Spencer map.}
\label{KSbacksection}

In the versal case, the tangent space to $K^g$ can be identified with $H^1(\mathcal X_0,\Theta_\pi)$ by use of the Kodaira-Spencer map. 

\begin{theorem}
\label{KSandCompletetheorem}
Let $\pi : \mathcal X_0\to\mathbb S^1$ be a proper CR submersion. Assume that for all $t\in\mathbb S^1$, the identity {\rm \eqref{versalcondition}} is fulfilled. Then there exists a fixed isomorphism
$$
\begin{CD}
H^1(\mathcal X_0,\Theta_\pi)@ >\varphi >> T_0K^g
\end{CD}
$$
such that the following property holds.

 Let $\Pi : \mathcal Z\to B$ be any deformation of $\mathcal X_0$. Let $\rho_0$ be its Kodaira-Spencer map at $0$ and let $f : (B,0)\to (K^g, 0)$ be given by completeness of $K^g$.  Then we have
$$
\varphi\circ \rho_0\equiv d_0f
$$
\end{theorem}

\begin{proof}
We have already seen that the natural projection
$$
\mathcal H^1=\bigcup_{t\in\mathbb S^1} H^1(X_t,\Theta_t)\longrightarrow \mathbb S^1
$$
can be endowed with a structure of a sheaf over the circle with stalk $H^1(X_t,\Theta_t)$ at $t$. With this structure, the set of smooth sections of $\mathcal H^1$ identifies with the cohomology group $H^1(\mathcal X_0,\Theta_\pi)$.\vspace{5pt}\\
The key point is that, when \eqref{versalcondition} is fulfilled, this sheaf is indeed a vector bundle over the circle and $H^1(\mathcal X_0,\Theta_\pi)$ identifies now with the set of smooth sections of this vector bundle, see \cite{K-S1}.\vspace{5pt}\\
If $\Pi : \mathcal Z\to B$ is a deformation of $\mathcal X_0$, recall that 
$$
\begin{CD}
\mathcal Y_t:=s^{-1}(t)@ >\Pi_t>>B
\end{CD}
$$
is a deformation of $X_t:=\pi^{-1}(t)\subset \mathcal X_0$, see \eqref{Yt}. Associated to it, when $B$ is finite-dimensional, we thus have a (classical) Kodaira-Spencer map
$$
\rho_t \ :\ T_0B\longrightarrow H^1(X_t,\Theta_t)
$$
and the family of these maps, when $t$ varies in $\mathbb S^1$, is exactly the Kodaira-Spencer map defined in \eqref{KSmap}.\vspace{5pt}\\
We want to apply these considerations to the Kuranishi family of $\mathcal X_0$. In this case, because of \eqref{CDbis}, the deformation $\mathcal Y_t$ reduces to a deformation over $(K_{c_0}, c_0(t))$. Indeed $\mathcal Y_t$ is equal to the pull-back of $\mathcal K_{c_0}\to (K_{c_0}, c_0(t))$ by the evaluation map at $t$. Hence $\rho_t$ gives a decomposition of the Kodaira-Spencer map $\rho$ of the family $\Pi^g : \mathcal K^g\to K^g$ into a family 
$$
\rho_t\ :\ T_{c_0(t)}K_{c_0}\longrightarrow H^1(X_t,\Theta_t)
$$
and we have a commutative diagram
$$
\begin{CD}
T_0K^g@ >\rho >> H^1(\mathcal X_0,\Theta_\pi)\cr
@ V ev_t VV @ VV ev_t V\cr
T_{c_0(t)}K_{c_0}@ >>\rho_t >H^1(X_t,\Theta_t)
\end{CD}
$$
Since we assume that \eqref{versalcondition} is fulfilled for all $t$, the space $K_{c_0}$ is versal at $c_0(t)$ and all the $\rho_t$ are isomorphisms. So is the map $\rho$.\vspace{5pt}\\
We define the map $\varphi$ of Theorem \ref{KSandCompletetheorem} to be $\rho^{-1}$.\vspace{5pt}\\
Now, because of \eqref{evaCD}, \eqref{CDbis} and recalling the markings, the second property of Theorem \ref{KSandCompletetheorem} is satisfied if and only if it is satisfied for each $t\in\mathbb S^1$, that is if
$$
\varphi_t\circ \rho_t\equiv d_0f_t
$$
which is true by the chain-property of the (classical) Kodaira-Spencer map.
\end{proof}

\section{Uniqueness of the Kuranishi type moduli space.}
\label{uniquesection}

In the classical setting, the versal property implies the uniqueness of the Kuranishi space as a germ. In our context, we get only uniqueness in the universal case. In the versal case, we have only a weaker form of uniqueness. Specifically, as a consequence of Theorem \ref{versaltheorem}, we have

\begin{corollary}
\label{uniquenesscorollary}
Let $\pi : \mathcal X_0\to\mathbb S^1$ be a proper CR submersion.
\begin{enumerate}
\item[(i)] Assume that the Kuranishi type family  $\mathcal K^g\to K^g$ is universal. Then it is unique, that is: if $\mathcal K^h\to K^h$ is another universal family for $\mathcal X_0$, then both families are isomorphic as germs.
\item[(ii)] Assume that $\mathcal K^g\to K^g$ is versal at $0$. Then it is unique in the following restricted sense: if $\mathcal K^h\to K^h$ is another versal family at $0$ such that $K^h$ is the loop space of the base $K'$ of some finite-dimensional analytic family $\mathcal K'\to K'$, then $K_{c_0}$ and $K'$ from the one hand, $K^g$ and $K^h$ from the other hand are isomorphic as germs.
\end{enumerate}
\end{corollary}

\begin{remark}
\label{Lindep}
In both cases, the Kuranishi type moduli space $K^g$ as well as $K_{c_0}$ are independent of the subspace $L$.
\end{remark}

\begin{remark}
In the versal case, comparing to the classical case of complex structures, one should expect uniqueness in the same general sense as in the universal case. However, the classical proofs of the implication: versality gives uniqueness use in a crucial way the finite-dimensionality of $K^g$ (see \cite[\S 1.6]{G-H-S}) and cannot be applied here.
\end{remark}

\begin{proof}

\noindent (i) This is just an application of the definition. Since both families are complete, there exist a morphism $f$ between the germ of family $\mathcal K^g$ and the germ of family $\mathcal K^h$, respectively $g$ between $\mathcal K^h$ and $\mathcal K^g$, both respecting the markings. Hence the 
germ of $\mathcal K^g$ is isomorphic to its pull back by $g\circ f$, and $\mathcal K^h$ to its pull back by $f\circ g$. By universality, both compositions must be the identity, hence $f$ and $g$ are isomorphisms.\medskip

\noindent (ii) Once again, since both families are complete, there exist a morphism $f$ between the germ of family $\mathcal K^g$ and the germ of family $\mathcal K^h$, respectively $g$ between $\mathcal K^h$ and $\mathcal K^g$, both respecting the markings. Theorem \ref{versaltheorem} implies that $K_{c_0}$ and $K'$ are versal at each point of $c_0$ (without loss of generality, we assume that both spaces contain $c_0$ and that both markings consist in the identification between $\mathcal X_0$ and $c_0$). From our definition of deformation and the particular form of $K^g$ and $K^h$, it follows that $f$ and $g$ are induced by analytic maps between $K_{c_0}$ and $K'$ that we still denote by $f$ and $g$. Hence, arguing as in \cite[\S 1.6]{G-H-S}, $f$ and $g$ are local isomorphisms at each point of $c_0$. Since they respect the markings, they are also bijective on $c_0$. Hence  the result.
\end{proof}

In the non-versal case, we lose uniqueness, that is both conclusions of Corollary \ref{uniquenesscorollary} may be false, see section \ref{exemples}.

\section{Rigidity and connectedness.}
\label{rigiditysection}
\subsection{Connectedness and extension of deformations}
\label{connectedss}
In the classical case of compact complex manifolds, Kuranishi's Theorem has as a consequence that every complex structure $J$ on $X^{diff}$ close enough to a fixed structure $J_0$ is connected to it. That is, there exists a 1-dimensional holomorphic (resp. smooth) deformation of $X_0$ that contains $X_J$. The proof just consists in choosing a disk (resp. a path) in the Kuranishi space of $X_0$ joining the base point to the point encoding $J$.\vspace{5pt}\\
In our case, the same result is true.

\begin{theorem}
\label{connectedtheorem}
Let $\pi : \mathcal X_0\to\mathbb S^1$ be a proper CR submersion, represented by an element $c_0$ of $\mathcal I(\mathcal X^{diff})$ satisfying Hypothesis {\rm \ref{noselfint}} as usual. 

Then, if $\mathcal X$ is a proper CR submersion close enough to $\mathcal X_0$ (that is, if $\mathcal X$ can be represented by a path $c$ close enough to $c_0$), there exists a holomorphic (resp. smooth) 1-dimensional deformation joining $\mathcal X_0$ and $\mathcal X$.
\end{theorem}

\begin{proof}
Since $\mathcal X$ is close to $\mathcal X_0$, it is represented by a point in the Kuranishi type moduli space $K^g$ of $\mathcal X_0$. That means that there exists a loop $c$ in $K_{c_0}$ encoding $\mathcal X$. Therefore, to construct a smooth deformation as desired, it is enough to construct an isotopy
$$
H\ :\ \mathbb S^1\times [0,1]\longrightarrow K_{c_0}
$$
such that
$$
H_0:=H(-,0)\equiv c_0\qquad\text{ and }\qquad H_1:=H(-,1)\equiv c.
$$
 Now, it is a classical fact that two smooth loops in an analytic space are isotopic {\it as soon as they are close enough one from the other}. Indeed, this is clear for complex manifolds. For analytic spaces, we may first desingularize and extend the loops we want to isotope. Observe the exceptional divisors being simply connected, there is no additional obstruction. Since we may assume $c$ to be arbitrarily close to $c_0$, the existence of $H$ follows.
\begin{remark}
There is no particular problem if the space $K_{c_0}$ is not reduced. Indeed, since we use maps from a reduced base (here a circle or an annulus) into $K_{c_0}$, they map into the reduction of $K_{c_0}$. Hence, we only have to desingularize the reduction of $K_{c_0}$. The same remark applies below and to the next results (Corollaries \ref{concorI} and \ref{concorII}).
\end{remark}
We treat the case of a $1$-dimensional {\it holomorphic} family joining $\mathcal X_0$ to $\mathcal X$. That amounts to finding some CR morphism of $\mathbb S^1\times\mathbb D$ in $K_{c_0}$ whose image contains $c_0$ and $c$. Once again, desingularizing and extending $c_0$ and $c$ if necessary, we may assume that $K_{c_0}$ is smooth. For each $\exp {2i\pi\theta}\in\mathbb S^1$, the corresponding holomorphic disk of $K_{c_0}$ must pass through $c_0(\exp {2i\pi\theta})$ and $c(\exp {2i\pi\theta})$. We can always construct such a disk $\mathbb D_\theta$ for $\theta$ fixed. And this can be done in a locally smooth way. The only problem that could appear is that, starting with $\mathbb D_0$ and constructing the family by extension, we finish with $\mathbb D_1$ different from $\mathbb D_0$. Lemma \ref{extensionlemma}, (ii) allows us to solve this problem.
\end{proof}

Indeed, we have even a stronger connectedness result.

\begin{corollary}
\label{concorI}
Let $\pi : \mathcal X_0\to\mathbb S^1$ be a proper CR submersion. 
If $X$ is a compact complex manifold close enough to some fiber $X_t$ of $\mathcal X_0$, then there exists a holomorphic (resp. smooth) $1$-dimensional deformation $\mathcal Z$ of $\mathcal X_0$ such that, for some $z$ in the base, the $t$-fiber of $\mathcal Z_z$ is biholomorphic to $X$.
\end{corollary}
\begin{proof}
Choose $c$ in $\mathcal I(\mathcal X^{diff})$ close to $c_0$ and satisfying that $X_{c(t)}$ is biholomorphic to $X$ and apply Theorem \ref{connectedtheorem}.
\end{proof}

Finally, we prove that a $1$-dimensional deformation of a fiber of $\mathcal X_0$ can be extended as a $1$-dimensional deformation of $\mathcal X_0$.

\begin{corollary}
\label{concorII}
Let $\pi : \mathcal X_0\to\mathbb S^1$ be a proper CR submersion. 
If $\mathcal Y\to B$ is a holomorphic (resp. smooth) deformation of some fiber $X_t$ of $\mathcal X_0$ over a $1$-dimensional reduced base, then there exists a holomorphic (resp. smooth) deformation $\mathcal Z$ of $\mathcal X_0$ over the same base inducing locally $\mathcal Y$, i.e. such that we have
$$
\restriction{(\mathcal Y_t)}{U}\equiv \mathcal Y
$$
for $U$ a neighborhood of $0$ in $B$.
\end{corollary}

Recall the definition \eqref{Yt} of $\mathcal Y_t$.
 
\begin{proof}
We just do the holomorphic case. Since $K_{c_0}$ is complete at $c_0(t)$, we may assume that $\mathcal Y$ is obtained by pull-back of $K_t$ (the Kuranishi space of $X_t$) over some disk of the base. As before, we may assume that $K_{c_0}$ is smooth. Let
$$
F\ :\ (\mathbb D,0)\subset (B,0)\longrightarrow (K_{c_0}, c_0(t))
$$
be the associated map. We just have to extend it into a CR map
$$
H\ :\ \mathbb S^1\times \mathbb D\longrightarrow K_{c_0}
$$
such that
$$
H(-,0)\equiv c_0\qquad\text{ and }\qquad H(t,-)\equiv F
$$
which is not really different from what we did in the proofs of Theorems \ref{versaltheorem} and Theorem \ref{connectedtheorem} thanks to Lemma \ref{extensionlemma}.
\end{proof}

\subsection{Rigidity}
\label{rigidityss}
As in the classical case, we say that 
\begin{definition}
\label{rigiditydef}
A proper CR submersion $\pi : \mathcal X_0\to\mathbb S^1$ is {\it rigid} if any deformation $\Pi : \mathcal Z\to B$ of $\mathcal X_0$ is (locally) isomorphic to a product
$\mathcal X_0\times B$.
\end{definition}
Here is a trivial example of a rigid CR submersion.

\begin{example}
\label{rigidexample}
Let $X_0$ be a rigid compact complex manifold (for example, $X$ is $\mathbb P^n$ for some $n>0$). Let $\mathcal X_0$ be a CR bundle with fiber $X_0$.
If this bundle is trivial, then the space $K_{c_0}$ is nothing else that a point. So is $K^g$. By Theorem \ref{main}, $\mathcal X_0$ is rigid.
If it is not trivial, then $K_{c_0}$ can be taken as the unit disk (cf. Theorem \ref{main}), but the family $\mathcal K_{c_0}\to K_{c_0}$ is trivial by construction and every $\mathcal X$ close to $\mathcal X_0$ is isomorphic to 
$$
X_0\times [0,1]/\sim \qquad\text{ where }\qquad (z,0)\sim (\phi (z),1)
$$
for $\phi$ a biholomorphism representing the monodromy of $\mathcal X_0$. Hence it is isomorphic to $\mathcal X_0$. The same argument shows that any deformation of $\mathcal X_0$ is trivial, hence $\mathcal X_0$ is rigid.
\end{example}

This is indeed the unique rigid example.

\begin{theorem}
\label{rigiditytheorem}
Let $\pi : \mathcal X_0\to\mathbb S^1$ be a proper CR submersion. Then $\mathcal X_0$ is rigid if and only if it is a CR bundle with rigid fiber.
\end{theorem}

\begin{proof}
We have already seen in the previous example that a CR bundle with rigid fiber is rigid. Conversely let $\mathcal X_0$ be rigid. If one of the fiber $X_t$ of $\mathcal X_0$ is not rigid, then by Corollary \ref{concorI}, there exists a non-trivial deformation of $\mathcal X_0$, so every fiber is rigid. By connectedness of the circle, this implies that all the fibers are biholomorphic, hence, by Fischer-Grauert's Theorem, it is a CR bundle. 
\end{proof}

\section{Examples}
\label{exemples}
\subsection{Trivial CR bundles}
\label{trivialCRexample}
Assume that $\pi : \mathcal X_0\to \mathbb S^1$ is a trivial CR bundle. Then $c_0$ is just a point of $\mathcal I$ and $K_{c_0}$ is the same as the Kuranishi space $K_0$ of the fiber $X_0$. So finally the Kuranishi type moduli space $K^g$ of $\mathcal X_0$ is $C^\infty (\mathbb S^1, K_0)$.
For example, if $\mathcal X_0$ is $\mathbb E_\tau\times\mathbb S^1$ (where $\mathbb E_\tau$ is the elliptic curve of modulus $\tau\in\mathbb H$), then $K_0$ is a neighborhood of $\tau$ in $\mathbb H$ and $K^g$ is the space of smooth maps from the circle to this neighborhood.

\subsection{Non-trivial CR bundles}
\label{badcaseexample}
Let $\omega$ be $\exp (2i\pi/3)$ and let $\mathbb E$ be the elliptic curve of modulus $\omega$.
Let $\mathcal X_0$ be the CR bundle with fiber $\mathbb E$ and monodromy $\omega$. Here $c_0$ is also a point, but we are in the special case of Theorem \ref{main} and we cannot take $K_{c_0}$ as a neighborhood of $\omega$ in $\mathbb H$. Let $\mathbb D_\omega$ be a disk centered at $\omega$ in $\mathbb H$. Let 
$$
V=\{z\in \mathbb C\mid \inf_{t\in [0,1]} \vert z-(t\omega+(1-t))\vert<\epsilon\}.
$$ 
This is a neighborhood of the segment joining $1$ to $\omega$ in $\mathbb C$. Let $\mathbb D_\epsilon$ be the open disk of radius $\epsilon$ centered at $1$. Observe that $\mathbb D_\epsilon$ is a neighborhood of $1$ included in $V$. Assume that $\epsilon$ is small enough to ensure that $\mathbb D_\epsilon$ and $\omega\mathbb D_\epsilon$ are disjoint.\vspace{5pt}\\
Following the proof of Theorem \ref{main}, we define $K_{c_0}$ to be
$$
K_{c_0}=\mathbb D_\omega\times V/\sim\qquad \text{ with } (\tau,w)\in \mathbb D_\omega\times\mathbb D_\epsilon\sim (\omega\cdot \tau, \omega w)
$$
where $\omega\cdot $ describes the action of the automorphism $\omega$ onto $\mathbb H$, that is
$$
\omega\cdot \tau=\dfrac{-1-\tau}{\tau}.
$$
So $K_{c_0}$ is biholomorphic to the product of a disk with an annulus. And $\mathcal X_0$ has non-trivial deformations, even if the situation may at first sight be rigid, due to the fact that no other elliptic curve than $\mathbb E$ admits $\omega$ as an automorphism. A CR submersion close to but different from $\mathcal X_0$ is encoded in a path $c$ in $\mathbb H$ with
$$
c(1)=\omega\cdot c(0).
$$
Such a structure is of course non-constant (in the sense that the fibers of the CR submersion are not all the same) and this explains how it is possible that the monodromy is not a biholomorphism of any fiber (cf. \cite{M-V}). This is indeed an example of a non-versal Kuranishi family.\vspace{5pt}\\
Observe that we are in the special case where the dimension of $K_{c_0}$ is one more than the dimension of the Kuranishi space it is constructed with. This extra-dimension comes from the fact that we need $\mathcal X_0$ to be represented by a path and not a loop in $K_{c_0}$. Without this trick, one should take as $K_{c_0}$ the quotient of $\mathbb D_\omega$ by the action generated by $\omega$; but we should then consider $K_{c_0}$ as {\it an orbifold}.\vspace{5pt}\\
More generally, if we take as $\mathcal X_0$ a CR bundle with fiber $X_0$ a compact complex $2$-torus and monodromy an automorphism of $X_0$ non-isotopic to the identity and non-periodic (such pairs exist, see \cite{G-V}), the same construction yields as $K_{c_0}$ the product of the Kuranishi space $K_0$ of $X_0$, an open set of $\mathbb C^4$, with an annulus. Forgetting the extra-dimension, one should take the quotient of $K_0$ by the action generated by this automorphism; but we should then consider a non-Hausdorff space. Once again, our Kuranishi type moduli space is not versal.\vspace{5pt}\\
On the contrary, when the monodromy of the fiber $X_0$ is not isotopic to the identity but extends as an automorphism of any manifold in the Kuranishi space $K_0$ of $X_0$, we may take $K_0$ as $K_{c_0}$ and gain one dimension with respect to our construction. However, with this "reduced" $K_{c_0}$, there is no equivalent to Theorem \ref{completetheorem}, because the space $\mathcal K^g$ is not complete. There is no hope to obtain a non-trivial CR bundle as pull-back by a {\it constant map}.\vspace{5pt}\\
Hence, we see that in all these examples, our space $K_{c_0}$ is not versal but it is minimal with respect to properties of Theorems \ref{main} and \ref{completetheorem}. This strongly suggests that there is no versal space in this situation.

\subsection{Hopf surfaces}
\label{Hopfexample}
Consider the quotient of 
$$
\mathbb C^2\setminus\{(0,0)\}\times \mathbb S^1
$$
by the action generated by the map
$$
(z,w,t)\longmapsto (2z+a(t)w^2, (2+b(t))w,t)
$$
 where $a$ and $b$ are two smooth functions from $\mathbb S^1\subset \mathbb C$ into $\mathbb R_{\geq 0}$ satisfying
\begin{equation}
\label{patha}
a(-1)=a(1)=0\qquad\text{ and }\qquad a(t)>0\text{ for }t\not =-1,1
\end{equation}
and
\begin{equation}
\label{pathb}
b(1)=0,\ b(-1)=2\qquad\text{ and }\qquad 0<b(t)<2\text{ for }t\not = -1,1.
\end{equation}
This defines a smoothly trivial CR submersion $\mathcal X_0$ over the circle with fibers (primary) Hopf surfaces. As usual, we denote by $X_t$ the fiber over $t$.\vspace{5pt}\\
The following facts are well-known (cf. \cite{Da} or \cite{We}).
\begin{enumerate}
\item[(i)] The Hopf surface $X_1$ is the quotient of $\mathbb C^2\setminus\{(0,0)\}$ by the  linear action generated by the matrix $2Id$. Its Kuranishi space is smooth of dimension four and a Kuranishi domain can be identified  with an open neighborhood $V$ of $2Id$ in $\text{GL}_2(\mathbb C)$ under the correspondence
$$
A\in V\longmapsto X_A:=\left (\mathbb C^2\setminus\{(0,0)\}\right )/\langle A\rangle .
$$
\item[(ii)] The Hopf surface $X_{-1}$ is the quotient of  the quotient of $\mathbb C^2\setminus\{(0,0)\}$ by the  linear action generated by the diagonal matrix with eigenvalues $2$ and $4$. Because $4$ is the square of $2$, we are in a resonant case and its Kuranishi space is a bit different from the previous one. It is smooth of dimension three and a Kuranishi domain can be identified with an open neighborhood $W$ of $(2,4,0)$ in $\mathbb C^3$ under the correspondence
$$
(\alpha,\beta,s)\in W\longmapsto X_{(\alpha,\beta,s)}:=\left (\mathbb C^2\setminus\{(0,0)\}\right )/\langle A_{(\alpha,\beta,s)}\rangle
$$
where
$$
A_{(\alpha,\beta,s)}(z,w)=(\alpha z+sw^2,\beta w).
$$
\item[(iii)] For $t\not = -1,1$, the Hopf surface $X_t$ is biholomorphic to the linear Hopf surface defined by the matrix
$$
A_\lambda (t)=
\begin{pmatrix}
2 &\lambda\cr
0 &2+b(t) 
\end{pmatrix}
$$
for any choice of $\lambda\in\mathbb C$. Its Kuranishi space is smooth of dimension two and a Kuranishi domain can be identified with an open neighborhood of $A_0(t)$ in the space of {\it diagonal} matrices.
\item[(iv)] The Kuranishi domain $V$ can be assumed to contain all matrices $A_{a(t)}(t)$ for $t$ different from $-1$. In the same way, the Kuranishi domain $W$ can be assumed to contain all triple $(2,2+b(t), a(t))$ for $t$ different from $1$.
\end{enumerate}
Taking these facts into account, we see that it is enough to choose adequately $V$ and $W$ and to glue $V$ to $W\times\mathbb C$ along a well-chosen open set to construct $K_{c_0}$.\vspace{5pt}\\
To be more precise, consider the closed path 
$$
t\in\mathbb S^1\cap \{\Re t\geq 0\}\longmapsto c^+(t):=A_{a(t)}(t)\in \text{GL}_2(\mathbb C).
$$
Let $\Vert -\Vert$ be the standard euclidian norm on $\mathbb C^4$. It induces a norm on $V\subset \text{GL}_2(\mathbb C)$ by identifying $\text{GL}_2(\mathbb C)$ to an open set of $\mathbb C^4$. Let $\epsilon$ be a positive real number and define
$$
V=\{A\in \text{GL}_2(\mathbb C)\mid \Vert A-c^+(t)\Vert <\epsilon\text{ for some }t\in\mathbb S^1\cap \{\Re t\geq 0\}\}.
$$
Symmetrically, consider the closed path
$$
t\in\mathbb S^1\cap \{\Re t\leq 0\}\longmapsto c^-(t):=(2, 2+b(t), a(t), 0)\in \mathbb C^4.
$$
Define 
$$
W=\{A\in\mathbb C^4\mid \Vert A-c^-(t)\Vert <\epsilon\text{ for some }t\in\mathbb S^1\cap \{\Re t\leq 0\}\}.
$$
Assume that $\epsilon$ is small enough so that
\begin{enumerate}
\item[(i)] The open sets $V$ and $W$ are Kuranishi domains.
\item[(ii)] The subset
$$
V_0=\{A\in \text{GL}_2(\mathbb C)\mid \Vert A-c^+(t)\Vert <\epsilon\text{ for }t=\pm i\}\subset V
$$
has two connected components and any point $A=(A_{ij})_{i,j=1,2}$ of $V_0$ satisfies
$$
0<\vert \lambda_1(A)\vert <\vert \lambda_2(A)\vert <4\quad\text{ and }\quad \vert A_{21}\vert >0
$$
where $\lambda_1(A)$ (respectively $\lambda_2(A)$) is the smallest (respectively biggest) eigenvalue of $A$.
\item[(iii)]  The subset
$$
W_0=\{A\in \mathbb C^4\mid \Vert A-c^-(t)\Vert <\epsilon\text{ for }t=\pm i\}\subset W
$$
has two connected component and any point $A=(A_i)_{i=1,\hdots, 4}$ of $W_0$ satisfies
$$
0<\vert A_1\vert <\vert A_2\vert <4\quad\text{ and }\quad \vert A_3\vert >0.
$$
\end{enumerate}
Then, because of the facts recalled above, the map
$$
A\in V_0\longmapsto (\lambda_1(A),\lambda_2(A), A_{21}, A_{12})\in W_0
$$
is an a.c. isomorphism. Gluing $V$ to $W$ through it gives the analytic space $K_{c_0}$ we are looking for. Here it is smooth of dimension four and it has the homotopy type of a circle. The paths \eqref{patha} and \eqref{pathb} glue together to give the path $c_0$.\vspace{5pt}\\
According to Theorem \ref{versaltheorem}, the associated space $K^g$ is not versal at $c_0$. This is easy to see in this case. Modifying the path \eqref{pathb} by replacing the zero fourth coordinate with any bump function (small enough in modulus) gives a path encoding $\mathcal X_0$ although it is different from $c_0$.\vspace{5pt}\\
Moreover, taking a bump function depending on a smooth parameter and performing the same construction, one obtains a trivial deformation of $\mathcal X_0$ over the interval with injective image in $K^g$. This contradicts versality.


\begin{thebibliography}{99}
\bibitem{B-S}
Berndtsson, B; Sibony, N.
\emph{The $\overline\partial$-equation on a positive current}.
Invent. Math. 147 (2002), 371--428.

\bibitem{Bog} Bogomolov, F.
\emph{Complex Manifolds and Algebraic Foliations}.
Publ. RIMS-1084, Kyoto, June 1996 (unpublished).

\bibitem{Bu}
Burel, T.
\emph{D\'eformations de feuilletages \`a feuilles complexes}.
Th\`ese de doctorat, Institut de Math\'ematiques de Bourgogne, Dijon, 2011.


\bibitem{Da}
Dabrowski, K.
\emph{Moduli spaces for Hopf surfaces}.
Math. Ann. 259 (1982), 201--225.

\bibitem{DG} Demailly, J.P.; Gaussier, H.
\emph{Algebraic embeddings of smooth almost complex structures}.
Preprint arxiv/1412.2899 (2014).

\bibitem{Do1}
Douady, A.
\emph{Le probl\`eme des modules pour les vari\'et\'es analytiques complexes}.
S\'em. Bourbaki 277 (1964/65).

\bibitem{Do2}
Douady, A.
\emph{Le probl\`eme des modules pour les sous-espaces analytiques compacts d'un espace analytique donn\'e}.
Ann. Inst. Fourier 16 (1966), 1--95.
 
\bibitem{EK-S}
El Kacimi, A.; Slim\`ene, J.
\emph{Cohomologie de Dolbeault le long des feuilles de certains feuilletages complexes}.
Ann. Inst. Fourier 60 (2010), 727--757.

\bibitem{F-K}
Forster, O.; Knorr, K.
\emph{Relativ-analytische R\"aume und die Koh\"arenz von Bildgarben}.
Invent. Math. 16 (1972), 113--160.

\bibitem{G-H-S}
Girbau, J.; Haefliger, A.; Sundararaman, D.
\emph{On deformations of transversely holomorphic foliations}.
J. Reine Angew. Math. 345 (1983), 122--147.

\bibitem{G-V}
Ghys, E; Verjovsky, A.
\emph{Locally free holomorphic actions of the complex affine group}.
Geometric study of foliations (Tokyo, 1993), World Sci. Publ., River Edge NJ, pp. 201--217, 1994.

\bibitem{Hamilton} Hamilton, R. S.
\emph{The inverse function theorem of Nash and Moser}.
 Bull. Amer. Math. Soc. (N.S.) 7 (1982), no. 1, 65--222.
 
 \bibitem{Hamilton2} Hamilton, R.S.
 \emph{Deformation Theory of Foliations}.
 Preprint, 88 pages, 1978 (unpublished).

\bibitem{H-I}
Henkin, I; Iordan, A.
\emph{Regularity of $\bar{\partial }$ on pseudoconcave compacts and applications}.
Asian J. Math. 4 (2000), 855--884.
\emph{Erratum}. 
Asian J. Math. 7 (2003), 147--148.

\bibitem{Ko}
Kodaira, K.
\emph{Complex manifolds and deformation of complex structures}.
Classics in Mathematics, Springer, Berlin, 2005.


\bibitem{K-S1}
Kodaira, K.; Spencer, D.C.
\emph{On deformations of complex analytic structures I}.
Ann. of Math. 67 (1958), 328--402.

\bibitem{K-S2}
Kodaira, K.; Spencer, D.C.
\emph{On deformations of complex analytic structures II}.
Ann. of Math. 67 (1958), 403--466.


\bibitem{Ku1}
Kuranishi, M.
\emph{On locally complete families of complex analytic structures}.
Ann. of Math. 75 (1962), 536--577.


\bibitem{Ku2}
Kuranishi, M.
\emph{New Proof for the Existence of Locally Complete Families of Complex Structures}.
Proc. Conf. Complex Analysis (Minneapolis, 1964), Springer, Berlin, pp. 142--154, 1965.

\bibitem{Ku3}
Kuranishi, M.
\emph{Deformations of Compact Complex Manifolds}.
Les presses de l'universit\'e de Montr\'eal, Montr\'eal,
1971.
\bibitem{Mathese}
Meersseman, L.
\emph{A new geometric construction of compact complex manifolds in any dimension}.
Math. Ann. 317 (2000), 79--115.

\bibitem{Me1}
Meersseman, L.
\emph{Foliated structure of the Kuranishi space and isomorphisms of deformation families of compact complex manifolds}.
Ann. Sci. \'Ec. Norm. Sup\'er. 44 (2011), 495--525.

\bibitem{Me2}
Meersseman, L.
\emph{Feuilletages par vari\'et\'es complexes et probl\`emes d'uniformisation}.
Panoramas \& Synth\`eses 34/35 (2011), 205--257.

\bibitem{M-V}
Meersseman, L.; Verjovsky, A.
\emph{On the moduli space of certain smooth codimension-one foliations of the $5$-sphere}.
J. Reine Angew. Math. 632 (2009), 143--202.

\bibitem{MK}
Morrow, J.; Kodaira, K.
\emph{Complex Manifolds}.
Holt, Rinehart and Winston, New York, 1971.

\bibitem{N-N}
Newlander, A; Nirenberg, L.
\emph{Complex analytic coordinates in almost complex manifolds}.
Ann. of Math. 65 (1957), 391--404.

\bibitem{Sc}
Schneider, M.
\emph{Halbstetigkeitss\"atze f\"ur relativ analytische R\"aume}.
Invent. Math. 16 (1972), 161--176.

\bibitem{Si}
Siu, Y.T.
\emph{Nonexistence of smooth Levi-flat hypersurfaces in complex projective spaces of dimension $\geqslant 3$}.
Ann. of Math. 151 (2000), 1217--1243.

\bibitem{Sl}
Slim\`ene, J.
\emph{Deux exemples de calcul explicite de cohomologie de Dolbeault feuillet\'ee}.
Proyecciones 27 (2008), 63--80.

\bibitem{Su}
Sundararaman, D.
\emph{Moduli, deformations and classifications of compact complex manifolds}.
Research Notes in Mathematics 45, Pitman, Boston, Mass.-London, 1980.

\bibitem{Wa1}
Wavrik, J.J.
\emph{Obstructions to the existence of a space of moduli}.
Global Analysis, Papers in honor of K. Kodaira, Princeton University Press, Princeton, pp. 403--414, 1969.

\bibitem{Wa2}
Wavrik, J.J.
\emph{Deforming cohomology classes}.
Trans. Amer. Math. Soc. 181 (1973), 341--350.

\bibitem{We}
Wehler, J.
\emph{Versal deformation of Hopf surfaces}.
J. Reine Angew. Math. 328 (1981), 22--32.
 \end{thebibliography}
\end{document}